\documentclass[a4paper,12pt]{article}
\usepackage{epsfig,amssymb,amsmath,amsthm,fullpage,wrapfig,pstricks}
\newtheorem{thm}{Theorem}

\newtheorem{cor}[thm]{Corollary}

\newtheorem{lem}[thm]{Lemma}
\newtheorem{propn}[thm]{Proposition}

\begin{document}

\title{Self-similarity and spectral asymptotics for the continuum
random tree}

\author{D.A. Croydon\footnote{Dept of Statistics,
University of Warwick, Coventry CV4 7AL, UK;
{d.a.croydon@warwick.ac.uk}.} and B. M. Hambly\footnote{Mathematical Institute, University of Oxford, 24-29 St Giles', Oxford, OX1 3LB, UK; hambly@maths.ox.ac.uk}}
\date{5 June 2007}
\maketitle

\begin{abstract}
We use the random self-similarity of the continuum random tree to show
that it is homeomorphic to a post-critically finite self-similar fractal
equipped with a random self-similar metric. As an application we
determine the mean and almost-sure leading order behaviour of the high
frequency asymptotics of the eigenvalue counting function associated
with the natural Dirichlet form on the continuum random tree. We also
obtain short time asymptotics for the trace of the heat semigroup and
the annealed on-diagonal heat kernel associated with this Dirichlet form.
\end{abstract}

\section{Introduction}\label{intro}

One of the reasons the continuum random tree of Aldous has attracted such
great interest is that it connects together a number of diverse areas of
probability theory. On one hand, it appears from discrete probability as the
scaling limit of combinatorial graph trees and probabilistic branching
processes; and on the other hand, it is intimately related with a continuous
time process, namely the normalised Brownian excursion, \cite{Aldous3}.
However, with both of these representations of the continuum random tree, there
does not appear to be an obvious description of the structure of the
set itself. In this paper we demonstrate that the continuum random
tree has a recursive description as a random self-similar fractal and
show that the set is always homeomorphic to a
deterministic subset of the Euclidean plane.
As an application of
this precise description of the random self-similarity of the continuum random tree, we deduce
results about the spectrum and on-diagonal heat kernel of the natural Dirichlet form on the set
using techniques developed for random recursive self-similar fractals.

From its graph tree scaling limit description, Aldous showed how the continuum
random tree has a certain random self-similarity, \cite{Aldous5}. In this article, we use
this result iteratively to label the continuum random tree,
$\mathcal{T}$, using a shift space over a three letter alphabet. This
enables us to show that there is an isometry from $\mathcal{T}$, with its natural metric
$d_{\mathcal{T}}$ (see Section \ref{crtdef} for a precise definition of
$\mathcal{T}$ and $d_\mathcal{T}$, and Section \ref{decompsec} for the decomposition of
$\mathcal{T}$ that we apply), to
a deterministic subset of
$\mathbb{R}^2$, $T$ say, equipped with a random metric $R$, $\mathbf{P}$-a.s.,
where $\mathbf{P}$ is the probability measure on the probability space upon
which all the random variables of the discussion are defined. This metric is
constructed using random scaling factors in an adaptation of the now
well-established techniques of \cite{Kigami} for building a resistance metric
on a post-critically finite self-similar fractal. We note that on a
tree the resistance and geodesic metrics are the same.
Furthermore, we show that
the isometry in question also links the natural Borel probability measures on
the spaces $(\mathcal{T},d_\mathcal{T})$ and $(T,R)$. The relevant measures will
be denoted by $\mu$ and $\mu^T$ respectively, with $\mu$ arising as the scaling
limit of the uniform measures on the graph approximations of $\mathcal{T}$ (see
\cite{Aldous3}, for example), and $\mu^T$ being the random self-similar measure
that is associated with the construction of $R$. The result that we prove is
the following; full descriptions of $(T,R,\mu^T)$ are given in Section \ref{selfsimsec},
and the isometry is defined in Section \ref{isosec}.

{\thm \label{first} There exists a deterministic post-critically finite
self-similar dendrite, $T$, equipped with a (random) self-similar metric, $R$,
and Borel probability measure, $\mu^T$, such that $(T,R,\mu^T)$ is equivalent
to $(\mathcal{T},d_\mathcal{T},\mu)$ as a measure-metric space,
$\mathbf{P}$-a.s.}
\bigskip


Previous analytic work on the continuum random tree in
\cite{Croydoncrt} obtained estimates on the quenched and the annealed
heat kernel for the tree. We can now adapt techniques of
\cite{Hamasymp} to consider the spectral asymptotics of the tree. As a
byproduct we are also able to refine the results on the annealed heat
kernel to show the existence of a short time limit for
$t^{2/3}\mathbf{E}p_t(\rho,\rho)$ at the root of the tree $\rho$, where the notation $\mathbf{E}$ is used to represent expectation under the probability measure $\mathbf{P}$.

The natural Dirichlet form on $L^2(\mathcal{T},\mu)$ may be thought of simply as
the electrical energy when we consider $(\mathcal{T},d_\mathcal{T})$ as a
resistance network. We shall denote this form by $\mathcal{E}_\mathcal{T}$, and
its domain $\mathcal{F}_\mathcal{T}$, and explain in Section \ref{crtdef} how it may be
constructed using results of \cite{Kigamidendrite}. The eigenvalues of
the triple $(\mathcal{E}_\mathcal{T},\mathcal{F}_\mathcal{T}, \mu)$
are defined to be the numbers $\lambda$ which satisfy
\begin{equation}\label{evaluedef}
\mathcal{E}_\mathcal{T}(u,v)=\lambda\int_\mathcal{T}uvd\mu,\hspace{20pt}\forall
v\in\mathcal{F}_\mathcal{T}
\end{equation}
for some eigenfunction
$u\in\mathcal{F}_\mathcal{T}$. The corresponding eigenvalue counting function,
$N$,  is obtained by setting
\begin{equation} \label{ecf}
N(\lambda):=\#\{\mbox{eigenvalues of
}(\mathcal{E}_\mathcal{T},\mathcal{F}_\mathcal{T}, \mu)\leq\lambda\},
\end{equation}
and we prove in Section \ref{specsec} that this is well-defined and finite
for any $\lambda\in\mathbb{R}$, $\mathbf{P}$-a.s. In Section \ref{specsec}, we also prove
the following result, which shows that asymptotically the mean and $\mathbf{P}$-a.s. behaviour
of $N$ are identical.

{\thm \label{second} There exists a deterministic constant $C_0\in(0,\infty)$ such that\\ (a)
$\lambda^{-2/3}\mathbf{E}N(\lambda)\rightarrow C_0$, as
$\lambda\rightarrow\infty$.\\ (b) $\lambda^{-2/3}N(\lambda)\rightarrow C_0$, as
$\lambda\rightarrow\infty$, $\mathbf{P}$-a.s.}
\bigskip 

To provide some context for this result, we will now briefly discuss some
related work. For the purposes of brevity, during the remainder of the
introduction, we shall use the notation $N(\lambda)$ to denote the eigenvalue
counting function of whichever problem is being considered. Classically, for
the usual Laplacian on a bounded domain $\Omega\subseteq\mathbb{R}^n$, Weyl's
famous theorem tells us that the eigenvalue counting function satisfies
\begin{equation}\label{weyl}
N(\lambda)=C_n|\Omega|\lambda^{n/2}+o(\lambda^{n/2}),\hspace{20pt}\mbox{as
}\lambda\rightarrow\infty, \end{equation} where $C_n$ is a constant depending
only on $n$, and $|\Omega|$ is the Lebesgue measure of $\Omega$, see \cite{Lap}. As a
consequence, in this setting, there exists a limit for
$\lambda^{-n/2}N(\lambda)$ as $\lambda\rightarrow\infty$. In the case
of deterministic p.c.f. self-similar fractals it is
known that
\[ N(\lambda)=\lambda^{d_S/2}(G(\ln\lambda)+o(1)),\hspace{20pt}\mbox{as
}\lambda\rightarrow\infty, \]
where $G$ is a periodic function, see
\cite{Kigami}, Theorem 4.1.5. The generic case has $G$ constant but
for fractals with a high degree of symmetry, such as the class of
nested fractals, (an example is the Sierpinski gasket), the function $G$
can be proved to be non-constant, and so no limit actually
exists for $\lambda^{-d_S/2}N(\lambda)$ as $\lambda\rightarrow\infty$
for these fractals. In the case of random recursive Sierpinski
gaskets, as studied in \cite{Hamasymp},
there are similar results, however the function $G$ must be multiplied
by a random weight variable, which can be thought of as a measure of
the volume of the fractal, and roughly corresponds to the factor
$|\Omega|$ in (\ref{weyl}). Again the generic case is that the limit of
the rescaled counting function exists and, in this setting, there are no known examples of
periodic behaviour. For the continuum random tree, no periodic fluctuations or random weight factors appear; this is due to
the non-lattice distribution of the Dirichlet
$(\frac{1}{2},\frac{1}{2},\frac{1}{2})$ random variables that are used
in the self-similar construction, and also the fact that summing the
three elements of the triple gives exactly one, $\mathbf{P}$-a.s.

It is also worth commenting upon the values of the exponent of $\lambda$ in the
leading order behaviour of $N(\lambda)$ in the classical and fractal setting.
From Weyl's result for bounded domains in $\mathbb{R}^n$, we see that the limit
\[d_S:=2\lim_{\lambda\rightarrow\infty}\frac{\ln N(\lambda)}{\ln\lambda}\]
is precisely $n$, matching the Hausdorff dimension of $\Omega$.
However, for deterministic and random self-similar fractals, this agreement is
not generally the case. For a large class of finitely ramified
fractals it has been proved that
\begin{equation}\label{spectraldim}
d_S=\frac{2d_H}{1+d_H},
\end{equation}
where $d_H$ is the Hausdorff dimension of the fractal in the
resistance metric (see \cite{Kigami}, Theorem 4.2.1, and \cite{Hamasymp},
Theorem 1.1). Due to its definition from the spectral asymptotics, the quantity $d_S$ has
become known as the spectral dimension of a (Laplacian on a) set. Clearly, from
the previous theorem, we see that for the continuum random tree $d_S=4/3$. This
result could have been predicted from the self-similar fractal picture of the
set given in Theorem \ref{first}, and (\ref{spectraldim}), noting that $d_H=2$ for the continuum random tree (see \cite{LegallDuquesne}). Observe that to be able to apply the result of \cite{LegallDuquesne}, the equivalence of the resistance and geodesic metrics on trees must be used.

Finally, let $X$ be the Markov process corresponding to
$(\mathcal{E}_\mathcal{T},\mathcal{F}_\mathcal{T}, \mu)$ and denote by
$(p_t(x,y))_{x,y\in\mathcal{T},\:t>0}$ its transition density;
 alternatively this is the heat kernel of the Laplacian associated with the Dirichlet
form. The existence of $p_t$ for $t>0$ was proved in \cite{Croydoncrt}, where
it was also shown that $t^{2/3}p_t(x,x)$ exhibits logarithmic fluctuations
globally, and log-logarithmic fluctuations for $\mu$-a.e. $x\in\mathcal{T}$, as
$t\rightarrow 0$. These fluctuations are caused by variations in the
``thickness'' of the measure $\mu$ over the space, which result in turn from
the randomness of the construction. However, the result of Theorem
\ref{second}(b) implies that these fluctuations must even out, when averaged
over the entire space. In particular, applying an Abelian theorem in the way
discussed in Remark 5.11 of \cite{Hamasymp}, we obtain the following limit
result for the trace of the heat semigroup, which we state without proof.

{\cor Let $C_0$ be the constant of Theorem \ref{second}, and $\Gamma$ be the
standard gamma function, then\\ (a)
\[t^{2/3}\mathbf{E}\int_\mathcal{T}p_t(x,x)\mu(dx)\rightarrow
C_0\Gamma(5/3),\hspace{20pt}\mbox{as }t\rightarrow 0,\] (b) $\mathbf{P}$-a.s.,
\[t^{2/3}\int_\mathcal{T}p_t(x,x)\mu(dx)\rightarrow
C_0\Gamma(5/3),\hspace{20pt}\mbox{as }t\rightarrow 0.\]}

Another corollary, which follows from the invariance under random
re-rooting of the continuum random tree, \cite{Aldous2}, allows us to
deduce from part (a) of this Corollary the following limit for the
annealed heat kernel at $\rho$,
the root of $\mathcal{T}$ (see Section \ref{crtdef} for a definition).
This tightens the result obtained in \cite{Croydoncrt}, Proposition
1.7 for the annealed heat kernel.

{\cor
 Let $C_0$ be the constant of Theorem \ref{second}, and $\Gamma$ be the
standard gamma function, then
\[ t^{2/3}\mathbf{E}p_t(\rho,\rho)\rightarrow C_0\Gamma(5/3) \mbox{ as }t\rightarrow
0, \]}

An outline of the paper is as follows. In Section 2 we introduce the
continuum random tree and give the natural Dirichlet form associated
with the tree. In Section 3 we use the decomposition of Aldous to give
a description of the tree via a sequence space. Once we have
established this we can map the continuum random tree into a
post-critically finite self-similar tree with a random metric. Finally
we show that the map ensures that the two sets are equivalent as
metric measure spaces. Once we have the picture as a self-similar set
with a random metric it is straightforward to deduce a decomposition
of the Dirichlet form and from this a natural scaling in the
eigenvalues. This leads to our results on the spectrum, and via
an Abelian theorem, to results on the trace of the heat semigroup.

\section{Continuum random tree}\label{crtdef}

The connection between trees and excursions is an area that has been
of much recent interest. In this section, we provide a brief
introduction to this link, a definition of the continuum random tree,
and also describe how to construct the natural Dirichlet form on this set.

We begin by defining the space of excursions, $U$, to be the set of
continuous functions $f:\mathbb{R}_+\rightarrow\mathbb{R}_+$ for which
there exists a $\tau(f)\in(0,\infty)$ such that $f(t)>0$ if and only
if $t\in(0,\tau(f))$. Given a function $f\in U$, we define a distance
on $[0,\tau(f)]$ by setting
\begin{equation}\label{distance}
d_f(s,t):=f(s)+f(t)-2m_f(s,t),
\end{equation}
where $m_f(s,t):=\inf\{f(r):\:r\in[s\wedge t,s\vee t]\}$. We then use the equivalence
\begin{equation}\label{eq}
s\sim t\hspace{20pt}\Leftrightarrow\hspace{20pt}d_f(s,t)=0,
\end{equation}
to define $\mathcal{T}_f:=[0,\tau(f)]/\sim$. Denoting by $[s]$ the
equivalence class containing $s$, it is elementary (see
\cite{LegallDuquesne}, Section 2) to check that
$d_{\mathcal{T}_f}([s],[t]):=d_f(s,t)$ defines a metric on
$\mathcal{T}_f$, and also that $\mathcal{T}_f$ is a {\it dendrite},
which is taken to mean a path-wise connected Hausdorff space
containing no subset homeomorphic to the circle. Furthermore, the
metric $d_{\mathcal{T}_f}$ is a shortest path metric on
$\mathcal{T}_f$, which means that it is additive along the paths of
$\mathcal{T}_f$. The {\it root} of the tree $\mathcal{T}_f$ is defined
to be the equivalence class $[0]$, and is denoted by $\rho_f$. A
natural volume measure to impose upon $\mathcal{T}_f$ is the
projection of Lebesgue measure on $[0,\tau(f)]$. In particular, for
open $A\subseteq\mathcal{T}_f$, let
\[\mu_f(A):=\ell\left(\{t\in[0,\tau(f)]:\:[t]\in A\}\right),\]
where, throughout this article, $\ell$ is the usual 1-dimensional
Lebesgue measure. This defines a Borel measure on
$(\mathcal{T}_f,d_{\mathcal{T}_f})$, with total mass equal to
$\tau(f)$.

We are now able to define the {\it continuum random
tree}\index{continuum random tree} as the random dendrite that we get
when the function $f$ is chosen according to the law of a suitably
scaled Brownian excursion. More precisely, we shall assume that there
exists an underlying probability space, with probability measure
$\mathbf{P}$, upon which is defined a process $W=(W_t)_{t=0}^1$ which
has the law of the normalised Brownian excursion, where, throughout
this article ``normalised'' is taken to mean ``scaled to return to the
origin for the first time at time 1''. In keeping with the notation
used so far in this section, the measure-metric space of interest
should be written $(\mathcal{T}_W, d_{\mathcal{T}_W}, \mu_W)$, the
distance on $[0,\tau(W)]$, defined at (\ref{distance}), $d_W$, and the
root, $\rho_W$. However, we shall omit the subscripts $W$ with the
understanding that we are discussing the continuum random tree in this
case. We note that $\tau(W)=1$, $\mathbf{P}$-a.s., and so
$[0,\tau(W)]=[0,1]$ and $\mu$ is a probability measure on
$\mathcal{T}$, $\mathbf{P}$-a.s. Moreover, that $\mu$ is non-atomic is
readily checked using simple path properties of $W$. Note that our
definition differs slightly from the Aldous continuum random tree,
which is based on the random function $2W$. Since this extra factor
only has the effect of increasing distances by a factor of 2, our
results are readily adapted to apply to Aldous' tree.

A further observation that will be useful to us is that between any
three points of a dendrite there is a unique branch-point. We shall
denote the branch-point of $x,y,z\in\mathcal{T}$ by $b(x,y,z)$, which
is the unique point in $\mathcal{T}$ lying on the arcs between $x$ and
$y$, $y$ and $z$, and $z$ and $x$.

Finally, we note that it is easy to check the conditions of
\cite{Kigamidendrite}, Theorem 5.4 to deduce that it is possible to
build a natural Dirichlet form on the continuum random tree.

{\thm
$\mathbf{P}$-a.s. there exists a local regular Dirichlet form
$(\mathcal{E}_\mathcal{T},\mathcal{F}_\mathcal{T})$ on
  $L^2(\mathcal{T},\mu)$, which is associated with the metric $d_\mathcal{T}$ through, for every $x\neq y$,
\begin{equation}\label{dtrecover}
d_\mathcal{T}(x,y)^{-1}=\inf\{\mathcal{E}_\mathcal{T}(f,f):\:f\in\mathcal{F}_\mathcal{T},\:f(x)=0,\:f(y)=1\}.
\end{equation}}
\bigskip

This final property means that the metric
$d_\mathcal{T}$ is indeed the resistance metric associated with
$(\mathcal{E}_\mathcal{T},\mathcal{F}_\mathcal{T})$. It will be the
eigenvalue counting function defined from
$(\mathcal{E}_\mathcal{T},\mathcal{F}_\mathcal{T},\mu)$ as at
(\ref{ecf}) for which we deduce asymptotic results in this article.

\section{Decomposition of the continuum random tree}\label{decompsec}

To make precise the decomposition of the continuum random tree that we
shall apply, we use the excursion description of the set introduced in
the previous section. This allows us to prove rigorously the
independence properties that are important to our argument. However,
it may not be immediately obvious exactly what the excursion picture
is telling us about the continuum random tree, and so, after Lemma
\ref{aldousdecomp}, we present a more heuristic discussion of the
procedure we use in terms of the related dendrites.

\begin{figure}[ht]
\centering
\includegraphics{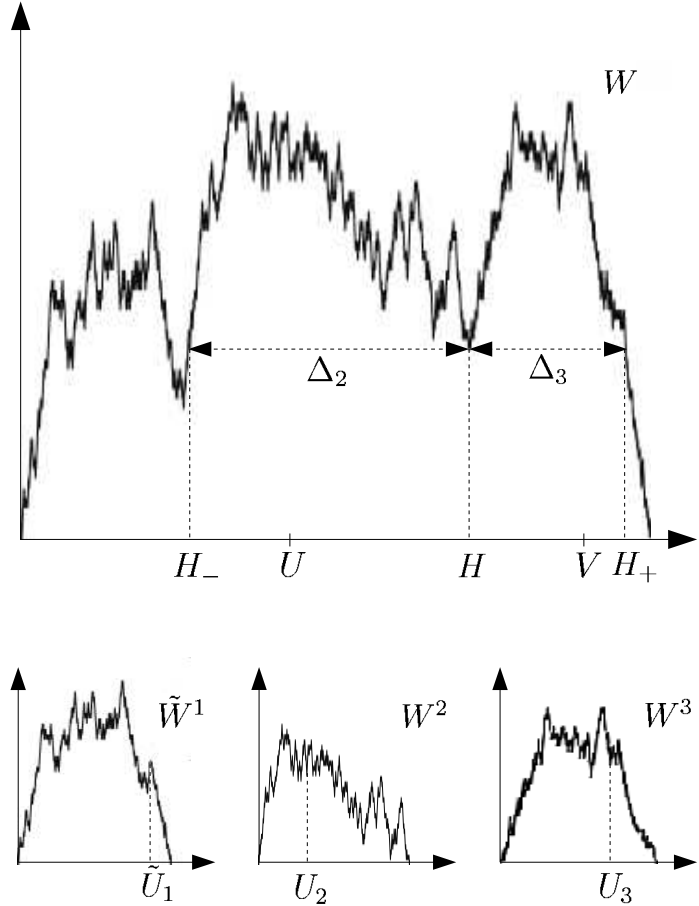}
\caption{Brownian excursion decomposition.}
\label{bed}
\end{figure}

The initial object of consideration is a triple
$(W,U,V)$, where $W$ is the normalised Brownian excursion, and $U$ and
$V$ are independent $U[0,1]$ random variables, independent of $W$. From this triple it is possible to
define three independent Brownian excursions. The following
decomposition is rather awkward to write down, but is made clearer by
Figure \ref{bed}. First, suppose $U<V$. On this set, it is
$\mathbf{P}$-a.s. possible to define $H\in [0,1]$ by
\begin{equation}\label{hdef}
\{H\}:=\{t\in[U,V]:\:W_t=\inf_{s\in[U,V]}W_s\}.
\end{equation}
We also define
\begin{equation}\label{hmindef}
H_-:=\sup\{t<U:\:W_t=W_H\},\hspace{20pt}H_+:=\inf\{t>V:\:W_t=W_H\},
\end{equation}
\[\Delta_1:=1+H_--H_+,\hspace{20pt}\Delta_2:=H-H_-,\hspace{20pt}\Delta_3:=H_+-H,\]
\[\tilde{U}_1:=\frac{H_-}{\Delta_1},\hspace{20pt}U_2:=\frac{U-H_-}{\Delta_2},\hspace{20pt}{U}_3:=\frac{V-H}{\Delta_3},\]
and for $t\in[0,1]$,
\[\tilde{W}^1_t:=\Delta_{1}^{-1/2}(W_{t\Delta_1}\mathbf{1}_{\{t\leq \tilde{U}_1\}}+W_{H_++(t-\tilde{U}_1)\Delta_1}\mathbf{1}_{\{t> \tilde{U}_1\}}),\]
\[{W}^2_t:=\Delta_{2}^{-1/2}(W_{H_-+t\Delta_2}-W_{H}),\]
\[{W}^3_t:=\Delta_{3}^{-1/2}(W_{H+t\Delta_3}-W_{H}).\]
Finally, it will be convenient to shift $\tilde{W}^1$ by $\tilde{U}_1$ so that the root of the corresponding tree is chosen differently. Thus, we define $W^1$ by
\[W^1_t:= \left\{\begin{array}{ll}W_{\tilde{U}_1}+W_{\tilde{U}_1+t}-2m(\tilde{U}_1,\tilde{U}_1+t),&\hspace{20pt}0\leq t \leq 1-\tilde{U}_1\\ W_{\tilde{U}_1}+W_{\tilde{U}_1+t-1}-2m(\tilde{U}_1+t-1,\tilde{U}_1),&\hspace{20pt}1-\tilde{U}_1\leq t \leq 1,\end{array}\right.\]
and set $U_1:=1-\tilde{U}_1$. If $U>V$, the definition of these quantities is similar, with $W^1$ again being the rescaled, shifted excursion containing $t=0$, $W^2$ being the rescaled excursion containing $t=U$, and $W^3$ being the rescaled excursion containing $t=V$. A minor adaptation of \cite{Aldous5}, Corollary 3, using the invariance under random re-rooting of the continuum random tree (see \cite{Aldous2}, Section 2.7), then gives us the following result, which we state without proof.

{\lem \label{aldousdecomp} The quantities $W^1,W^2,W^3,U_1,U_2,U_3$ and $(\Delta_1,\Delta_2,\Delta_3)$ are independent. Each $W^i$ is a normalised Brownian excursion, each $U_i$ is $U[0,1]$, and $(\Delta_1,\Delta_2,\Delta_3)$ has the Dirichlet $(\frac{1}{2},\frac{1}{2},\frac{1}{2})$ distribution.}
\bigskip 

Describing the result in terms of the corresponding trees gives a much
clearer picture of what the above decomposition does. Using the
notation of Section \ref{crtdef}, let $(\mathcal{T},d_\mathcal{T},
\mu)$ be the continuum random tree associated with $W$, and $\rho=[0]$
its root. Again, we use $[t]$, for $t\in[0,1]$, to represent the
equivalence classes of $[0,1]$ under the equivalence relation defined
at (\ref{eq}). If we define $Z^1:=[U]$ and $Z^2:=[V]$, then $Z^1$ and
$Z^2$ are two independent $\mu$-random vertices of $\mathcal{T}$. We
now split the tree $\mathcal{T}$ at the branch-point
$b(\rho,Z^1,Z^2)$, which may be checked to be equal to $[H]$, and
denote by $\mathcal{T}^1$, $\mathcal{T}^2$ and $\mathcal{T}^3$ the
components of $\mathcal{T}$ containing $\rho$, $Z^1$ and $Z^2$
respectively. Choose the root of each subtree to be equal to
$b(\rho,Z^1,Z^2)$ and, for $i=1,2,3$, let $\mu^i$ be the probability
measure on $\mathcal{T}^i$ defined by $\mu^i(A)=\mu(A)/\Delta_i$, for
measurable $A\subseteq\mathcal{T}^i$, where
$\Delta_i:=\mu(\mathcal{T}^i)$. The previous result tells us precisely
that $(\mathcal{T}^i,\Delta_i^{-1/2}d_\mathcal{T},\mu^i)$, $i=1,2,3$,
are three independent copies of
$(\mathcal{T},d_\mathcal{T},\mu)$. Furthermore, if $Z_i:=\rho,Z^1,Z^2$
for $i=1,2,3$, respectively, then $Z_i$ is a $\mu^i$-random variable
in $\mathcal{T}^i$. Finally, all these quantities are independent of
the masses
$(\mu(\mathcal{T}^1),\mu(\mathcal{T}^2),\mu(\mathcal{T}^3))$, which
form a Dirichlet $(\frac{1}{2},\frac{1}{2},\frac{1}{2})$
triple. Although it is possible to deal with the subtrees directly
using conditional definitions of the random variables to decompose the
continuum random tree in this way, the excursion description allows us
to keep track of exactly what is independent more easily, and it is to
this setting that we return. However, we shall not completely neglect
the tree description of the algorithm we now introduce, and a summary
in this vein appears after Proposition \ref{decompprop}.

We continue by applying inductively the decomposition map from
$U^{(1)}\times[0,1]^2$ to ${U^{(1)}}^3\times [0,1]^3\times\Delta$
(where $\Delta$ is the standard 2-simplex) that takes the triple
$(W,U,V)$ to the collection
$(W^1,W^2,W^3,U_1,U_2,U_3,(\Delta_1,\Delta_2,\Delta_3))$ of excursions
and uniform and Dirichlet random variables. We shall denote this
decomposition map by $\Upsilon$. To label objects in our consideration
it will be useful to use, as an address space, sequences of
$\{1,2,3\}$. In particular, we will write the collections of finite
sequences as, for $n\geq 0$,
\[\Sigma_n:=\{1,2,3\}^n,\hspace{20pt}\Sigma_*:=\bigcup_{m\geq
0}\Sigma_m,\] where $\Sigma_0:=\{\emptyset\}$. Later, we will refer to
the space of infinite sequences of $\{1,2,3\}$, which we denote by
$\Sigma$, and also apply some further notation, which we introduce
now. For $i\in\Sigma_m, j\in\Sigma_n, k\in\Sigma$, write $ij=i_1\dots
i_m j_1 \dots j_n$, and $ik=i_1\dots i_m k_1 k_2 \dots$. For
$i\in\Sigma_*$, denote by $|i|$ the integer $n$ such that
$i\in\Sigma_n$ and call this the {\it length}\index{length} of
$i$. For $i\in \Sigma_n\cup\Sigma$, $n\geq m$, the {\it
truncation}\index{truncation} of $i$ to length $m$ is written as
$i|m:=i_1\dots i_m$.

Now, suppose we are given an independent collection
$(W,U,(V_i)_{i\in\Sigma_*})$, where $W$ is a normalised Brownian
excursion, $U$ is $U[0,1]$, and $(V_i)_{i\in\Sigma_*}$ is a family of
independent $U[0,1]$ random variables. Set
$(W^\emptyset,U_\emptyset):=(W,U)$. Given $(W^i, U_i)$, define
\[(W^{i1},W^{i2},W^{i3},U_{i1},U_{i2},U_{i3},(\Delta_{i1},\Delta_{i2},\Delta_{i3})):=\Upsilon(W^i,U_i,V_i),\]
and denote the filtration associated with
$(\Delta_i)_{i\in\Sigma_*\backslash\{\emptyset\}}$ by
$(\mathcal{F}_n)_{n\geq0}$. In particular,
$\mathcal{F}_n:=\sigma(\Delta_i:\:|i|\leq n)$. The subsequent result
is easily deduced by applying the previous lemma repeatedly.

{\thm For each $n$, $((W^i,U_i, V_i))_{i\in\Sigma_n}$ is an independent collection of independent triples consisting of a normalised Brownian excursion and two $U[0,1]$ random variables, and moreover, the entire family of random variables is independent of $\mathcal{F}_n$.}
\bigskip 

Resulting from this construction, the collection $(\Delta_i)_{i\in\Sigma_*\backslash\{\emptyset\}}$ has some
particularly useful independence properties, which we will use in the
next section to build a random self-similar fractal related to
$\mathcal{T}$. Furthermore, Lemma \ref{aldousdecomp} implies that each
triple of the form $(\Delta_{i1},\Delta_{i2},\Delta_{i3})$ has the
Dirichlet $(\frac{1}{2},\frac{1}{2},\frac{1}{2})$
distribution. Subsequently we will also be interested in the
collection $(w(i))_{i\in\Sigma_*\backslash\{\emptyset\}}$, where for
each $i$, we define
\[w(i):=\Delta_i^{1/2},\] and will write $l(i)$ to represent the
product $w(i|1)w(i|2)\dots w(i||i|)$, where $l(\emptyset):=1$. The
reason for considering such families is that, in our decomposition of
the continuum random tree,
$(\Delta_i)_{i\in\Sigma_*\backslash\{\emptyset\}}$ and
$(w(i))_{i\in\Sigma_*\backslash\{\emptyset\}}$ represent the mass and
length scaling factors respectively.

By viewing the inductive procedure for decomposing excursions as the
repeated splitting of trees in the way described after Lemma
\ref{aldousdecomp}, it is possible to use the above algorithm to break
the continuum random tree into smaller components, with the subtrees
in the $n$th level of construction being described by the excursions
$(W^i)_{i\in\Sigma_n}$. The maps we now introduce will make this idea
precise. For the remainder of this section, the arguments that we give
hold $\mathbf{P}$-a.s. First, denote by $H^i$, $H_-^i$ and $H_+^i$ the
random variables in $[0,1]$ associated with $(W^i,U_i,V_i)$ by the
formulae at (\ref{hdef}) and (\ref{hmindef}). Let
$i\in\Sigma_*$. Define, for $t\in[0,1]$,
\[\phi_{i1}(t):=(H_+^i+t\Delta_{i1})\mathbf{1}_{\{t<U_{i1}\}}+(t-U_{i1})\Delta_{i1}\mathbf{1}_{\{t\geq
U_{i1}\}},\] and if $U_i<V_i$, define $\phi_{i2}$ and $\phi_{i3}$ to
be the linear contractions from $[0,1]$ to $[H_-^i,H^i]$ and
$[H^i,H_+^i]$ respectively. If $U_i>V_i$, the images of $\phi_{i2}$
and $\phi_{i3}$ are reversed. Note that, for each $i$, the map
$\phi_{i}$ satisfies, for any measurable $A\subseteq[0,1]$,
\begin{equation}\label{measurescal}
\ell(\phi_i(A))=\Delta_i\ell(A),
\end{equation}
where $\ell$ is the usual Lebesgue measure on $[0,1]$. Importantly,
these maps also satisfy a certain distance scaling property. In
particular, it is elementary to check from the definitions of the
excursions that, for any $i\in\Sigma_*$, $j\in \{1,2,3\}$,
\begin{equation}\label{disscal}
d_{W^i}(\phi_{ij}(s),\phi_{ij}(t))=w(ij)d_{W^{ij}}(s,t),\hspace{20pt}\forall s,t\in [0,1],
\end{equation}
where $d_{W^i}$ is the distance on $[0,1]$ associated with $W^i$ by
the definition at (\ref{distance}). This equality allows us to define
a map on the trees related to the excursions. Let
$(\tilde{\mathcal{T}}_i,d_{\tilde{\mathcal{T}}_i})$ be the metric
space dendrite determined from $W^i$ by the equivalence relation given
at (\ref{eq}). Denote the corresponding equivalence classes $[t]_i$
for $t\in[0,1]$. Now define, for $i\in\Sigma_*$, $j\in \{1,2,3\}$,
\begin{eqnarray*}
\tilde{\phi}_{ij}:\tilde{\mathcal{T}}_{ij}&\rightarrow&\tilde{\mathcal{T}}_i\\
\mbox{$[t]_{ij}$}&\mapsto&[\phi_{ij}(t)]_{i}.
\end{eqnarray*}
The following result is readily deduced from the distance scaling
property at (\ref{disscal}), and so we state it without proof.

{\lem\label{scalyo} $\mathbf{P}$-a.s., for every $i\in\Sigma_*$, $j\in \{1,2,3\}$, $\tilde{\phi}_{ij}$ is well-defined and moreover,
\[d_{\tilde{\mathcal{T}}_i}(\tilde{\phi}_{ij}(x),\tilde{\phi}_{ij}(y))=w(ij)d_{\tilde{\mathcal{T}}_{ij}}(x,y),\hspace{20pt}\forall
  x,y\in \tilde{\mathcal{T}}_{ij}.\]}

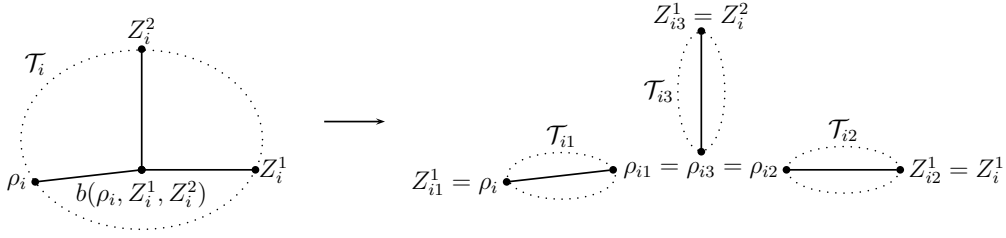
\begin{figure}[t]
\centering
\scalebox{.8}{
\begin{pspicture}(17,4)
\psellipse[linestyle=dotted](2.5,1.5)(2,1.5)
\psline[arrows=*-*](2.5,1)(4.37,1)
\psline[arrows=*-*](2.5,1)(2.5,3)
\psline[arrows=*-*](2.5,1)(.75,.8)
\rput(2.5,.6){$b(\rho_i,Z_{i}^1,Z_i^2)$}
\rput(4.67,1){$Z_{i}^1$}
\rput(2.5,3.3){$Z_i^2$}
\rput(.45,.8){$\rho_i$}
\rput(.75,2.75){$\mathcal{T}_i$}
\psline[arrows=->](5.5,1.8)(6.5,1.8)
\psline[arrows=*-*](8.5,.8)(10.25,1)
\rput(11.7,1){$\rho_{i1}=\rho_{i3}=\rho_{i2}$}
\psline[arrows=*-*](13.1,1)(14.97,1)
\psline[arrows=*-*](11.7,1.3)(11.7,3.3)
\rput(7.65,.8){$Z_{i1}^1=\rho_i$}
\rput(15.9,1){$Z_{i2}^1=Z_i^1$}
\rput(11.7,3.55){$Z_{i3}^1=Z_i^2$}
\psellipse[linestyle=dotted](14.035,1)(.935,.4)
\psellipse[linestyle=dotted](11.7,2.3)(.4,1)
\psellipse[linestyle=dotted](9.4,.9)(.9,.4)
\rput(9.4,1.55){$\mathcal{T}_{i1}$}
\rput(14.035,1.65){$\mathcal{T}_{i2}$}
\rput(11,2.3){$\mathcal{T}_{i3}$}
\end{pspicture}}
\caption{Continuum random tree decomposition.}
\label{crtd}
\end{figure}

By iterating the functions
$(\tilde{\phi}_{i})_{i\in\Sigma_*\backslash\{\emptyset\}}$, we can map
any $\tilde{\mathcal{T}}_i$ to the original continuum random tree,
$\mathcal{T}\equiv\tilde{\mathcal{T}}_\emptyset$, which is the object
of interest. We will denote the map from $\tilde{\mathcal{T}}_i$ to
$\mathcal{T}$ by
$\tilde{\phi}_{*i}:=\tilde{\phi}_{i|1}\circ\tilde{\phi}_{i|2}\circ\dots\circ\tilde{\phi}_{i}$,
and its image by
$\mathcal{T}_i:=\tilde{\phi}_{*i}(\tilde{\mathcal{T}}_i)$. It is these
sets that form the basis of our decomposition of $\mathcal{T}$. We
will also have cause to refer to the following points in
$\mathcal{T}_i$:
\[\rho_i:=\tilde{\phi}_{*i}([0]_i),\hspace{20pt}Z^1_i:=\tilde{\phi}_{*i}([U_i]_i),\hspace{20pt}Z^2_i:=\tilde{\phi}_{*i}([V_i]_i).\]
Although it has been quite hard work arriving at the definition of
$(\mathcal{T}_i)_{i\in\Sigma_*}$, the properties of this family of
sets that we will need are derived without too many difficulties from
the construction. The proposition we now prove includes the following
results: the sets $(\mathcal{T}_i)_{i\in\Sigma_n}$ cover
$\mathcal{T}$; $\mathcal{T}_i$ is simply a rescaled copy of
$\tilde{\mathcal{T}}_i$ with $\mu$-measure $l(i)^2$; the overlaps of
sets in the collection $(\mathcal{T}_i)_{i\in\Sigma_n}$ are small; and
also describes various relationships between points of the form
$\rho_{i}$, $Z_i^1$ and $Z_i^2$. This result is summarised in Figure \ref{crtd}.

{\propn \label{decompprop} $\mathbf{P}$-a.s., for every $i\in\Sigma_*$,\\
(a) $\mathcal{T}_i=\cup_{j\in\Sigma_n}\mathcal{T}_{ij}$, for all $n\geq 0$.\\
(b) $(\mathcal{T}_i,d_{\mathcal{T}})$ and
  $(\tilde{\mathcal{T}}_i,l(i)d_{\tilde{\mathcal{T}}_i})$ are
  isometric.\\
(c) $\rho_{i1}=\rho_{i2}=\rho_{i3}=b(\rho_{i},Z_i^1,Z_i^2)$.\\
(d) $Z_{ij}^1=\rho_{i},Z_i^1,Z_i^2$, for $j=1,2,3$ respectively.\\
(e) $\rho_{i}\not\in\mathcal{T}_{i2}\cup\mathcal{T}_{i3}$, $Z^1_{i}\not\in\mathcal{T}_{i1}\cup\mathcal{T}_{i3}$ and $Z^2_{i}\not\in\mathcal{T}_{i1}\cup\mathcal{T}_{i2}$.\\
(f) if $|j|=|i|$, but $j\neq i$, then $\mathcal{T}_i\cap\mathcal{T}_j=\{\rho_i\}$ when $j|(|j|-1)=i|(|i|-1)$, and $\mathcal{T}_i\cap\mathcal{T}_j=\emptyset$ or $\{Z_i^1\}$ otherwise.\\
(g) $\mu(\mathcal{T}_i)=l(i)^2$.}
\begin{proof} By induction, it suffices to show that (a) holds for $n=1$. By definition, we have $\cup_{j\in \{1,2,3\}}\phi_{ij}([0,1])=[0,1)$, and so
\[\tilde{\mathcal{T}}_i=\cup_{j\in
\{1,2,3\}}\{[\phi_{ij}(t)]_i:\:t\in[0,1]\}=\cup_{j\in
\{1,2,3\}}\tilde{\phi}_{ij}(\tilde{\mathcal{T}}_{ij}),\] where we
apply the definition of $\tilde{\phi}_{ij}$ for the final
equality. Applying $\tilde{\phi}_{*i}$ to both sides of this equation
completes the proof of (a). Part (b) is an immediate consequence of
the definition of $\mathcal{T}_i$ and the distance scaling property of
$\tilde{\phi}_{*i}$ proved in Lemma \ref{scalyo}.

Analogous to the remark made after Lemma \ref{aldousdecomp}, the point $[H^i]_i$ represents the branch-point of $[0]_i$, $[U_i]_i$ and $[V_i]_{i}$ in $\tilde{\mathcal{T}}_i$. Thus, since $\tilde{\phi}_{*i}$ is simply a rescaling map, we have that
\[b(\rho_{i},Z_i^1,Z_i^2)=b(\tilde{\phi}_{*i}([0]_i),\tilde{\phi}_{*i}([U_i]_i),\tilde{\phi}_{*i}([V_i]_i))=\tilde{\phi}_{*i}([H^i]_i).\]
Now, note that for any $j\in \{1,2,3\}$, we have by definition that
$\phi_{ij}(0)\in\{H^i,H^i_-,H^i_+\}$, and so
$[\phi_{ij}(0)]_i=[H^i]_i$. Consequently,
\begin{equation}\label{groggy}
\tilde{\phi}_{*i}([H^i]_i)=\tilde{\phi}_{*i}([\phi_{ij}(0)]_i)=\tilde{\phi}_{*ij}([0]_{ij})=\rho_{ij},
\end{equation}
which proves (c). Part (d) and (e) are easy to check from the construction using similar ideas and so their proof is omitted.

Now note that, for $k\in\Sigma_*$, the decomposition of the
excursions, and the fact that the local minima of a Brownian excursion
are distinct, implies that for $j_1,j_2\in \{1,2,3\}$, $j_1\neq j_2$,
we have
$\tilde{\phi}_{kj_1}(\tilde{\mathcal{T}}_{kj_1})\cap\tilde{\phi}_{kj_2}(\tilde{\mathcal{T}}_{kj_2})=\{[H^k]_k\}$. Applying
the injection $\tilde{\phi}_{*k}$ to this equation yields
\begin{equation}\label{intersection}
\mathcal{T}_{kj_1}\cap\mathcal{T}_{kj_2}=\{\tilde{\phi}_{*k}([H^k]_k)\}=\{\rho_{k1}\},
\end{equation}
with the second equality following from (\ref{groggy}). This fact will
allow us to prove (f) by induction on the length of $i$. Obviously,
there is nothing to prove for $|i|=0$. Suppose now that $|i|\geq 1$
and the desired result holds for any index of length strictly less
than $|i|$. Suppose $|j|=|i|$, but $j\neq i$, and define
$k:=i|(|i|-1)$. If $j|(|j|-1)\neq k$, then the inductive hypothesis
implies that
$\mathcal{T}_i\cap\mathcal{T}_j\subseteq\mathcal{T}_k\cap\mathcal{T}_{j|(|j|-1)}\subseteq\{\rho_k,Z_k^1\}$,
where we apply part (a) to obtain the first inclusion. Using parts (d)
and (e) of the proposition it is straightforward to deduce from this
that $\mathcal{T}_i\cap\mathcal{T}_j\subseteq\{Z_i^1\}$ in this
case. If $j|(|j|-1)=k$, then we can apply the equality at
(\ref{intersection}) to obtain that
$\mathcal{T}_i\cap\mathcal{T}_j=\{\rho_{k1}\}=\{\rho_i\}$, which
completes the proof of part (f).

Finally, $\mu$ is non-atomic and so
$\mu(\mathcal{T}_i)=\mu(\mathcal{T}_i\backslash\{\rho_i,Z_i^1\})$. Hence,
by the disjointness of the sets and the fact that $\mu$ is a
probability measure, we have
$1\geq\sum_{i\in\Sigma_n}\mu(\mathcal{T}_i\backslash\{\rho_i,Z_i^1\})=\sum_{i\in\Sigma_n}\mu(\mathcal{T}_i)$. Now,
by definition, for each $i$,
\[\mathcal{T}_i=\{\tilde{\phi}_{*i}([t]_i):\:t\in[0,1]\}=\{[t]:\:t\in\phi_{i|1}\circ\phi_{i|2}\circ\dots\circ\phi_{i}([0,1])\}.\]
Thus, since $\mu$ is the projection of Lebesgue measure, this implies
that $\mu(\mathcal{T}_i)$ is no smaller than
$\ell(\phi_{i|1}\circ\phi_{i|2}\circ\dots\circ\phi_{i}([0,1]))$. By
repeated application of (\ref{measurescal}), this lower bound is equal
to $\Delta_{i|1}\Delta_{i|2}\dots\Delta_i=l(i)^2$. Now observe that,
because $(\Delta_{i1},\Delta_{i2},\Delta_{i3})$ are Dirichlet
$(\frac{1}{2},\frac{1}{2},\frac{1}{2})$ random variables, we have
$\Delta_{i1}+\Delta_{i2}+\Delta_{i3}=1$ for every $i\in\Sigma_*$, and
from this it is simple to show that
$\sum_{i\in\Sigma_n}l(i)^2=1$. Hence
$\sum_{i\in\Sigma_n}\mu(\mathcal{T}_i)\geq\sum_{i\in\Sigma_n}l(i)^2=1$. Thus
$\sum_{i\in\Sigma_n}\mu(\mathcal{T}_i)$ is actually equal to 1, and
moreover, (g) must hold.
\end{proof}

With regards to Figure \ref{crtd}, note that the fact
that sets from $(\mathcal{T}_{ij})_{j\in \{1,2,3\}}$ only intersect at
$\rho_{i1}$ is proved in part (f) of the above proposition, and so the diagram
is representative of the set structure of the
decomposition. Furthermore, it is clear that the sets $\mathcal{T}_i$
are all compact dendrites, because they are simply rescaled versions
of the compact dendrites $\tilde{\mathcal{T}}_i$.

The tree description of the inductive algorithm runs as
follows. Suppose that the triples
$((\mathcal{T}_i,l(i)^{-1}d_{\mathcal{T}},\mu^i))_{i\in\Sigma_n}$ are
independent copies of $(\mathcal{T},d_{\mathcal{T}},\mu)$, independent
of $\mathcal{F}_n$, where $\mu^i(A):=\mu(A)/\mu(\mathcal{T}_i)$ for
measurable $A\subseteq\mathcal{T}_i$. Furthermore, suppose
$\mathcal{T}_i$ has root $\rho_i$, and $Z^1_i$ and $Z^2_i$ are two
$\mu^i$-random variables in $\mathcal{T}_i$. For $j=1,2,3$, define
$\mathcal{T}_{ij}$ to be the component of $\mathcal{T}_i$ (when split
at $b(\rho_i,Z_i^1,Z_i^2)$) containing $\rho_i,Z^1_i,Z^2_i$
respectively. Define $\Delta_{ij}:=\mu^i(\mathcal{T}_{ij})$, and equip
the sets with the metrics
$\Delta_{ij}^{-1/2}l(i)^{-1}d_{\mathcal{T}}=l(ij)^{-1}d_{\mathcal{T}}$
and measures $\mu^{ij}$, defined by
\[\mu^{ij}(A):=\frac{\mu^i(A)}{\Delta_{ij}}=\frac{\mu(A)}{\mu(\mathcal{T}_{ij})}.\]
Then the triples
$((\mathcal{T}_i,l(i)^{-1}d_{\mathcal{T}},\mu^i))_{i\in\Sigma_{n+1}}$
are independent copies of the continuum random tree, independent of
$\mathcal{F}_{n+1}$. Moreover, for $i\in\Sigma_{n+1}$, the algorithm
gives us the root $\rho_i$ of $\mathcal{T}_i$ and also a
$\mu^i$-random vertex, $Z_i^1$. To continue the algorithm, we pick
independently for each $i\in\Sigma_{n+1}$ a second $\mu^i$-random
vertex, $Z_i^2$. Note that picking this extra $\mu^i$-random vertex is
the equivalent of picking the $U[0,1]$ random variable $V_i$ in the
excursion picture.

To complete this section, we introduce one further family of variables
associated with the decomposition of the continuum random tree. From
Proposition \ref{decompprop}(f), observe that the sets in
$(\mathcal{T}_i)_{i\in\Sigma_n}$ only intersect at points of the form
$\rho_i$ or $Z_i^1$. Consequently, it is possible to consider the two
point set $\{\rho_i,Z_i^1\}$ to be the boundary of
$\mathcal{T}_i$. Denote the renormalised distance between boundary
points by, for $i\in \Sigma_*$,
\[D_i:=l(i)^{-1}d_\mathcal{T}(\rho_i,Z_i^1).\] By construction, we
have that $d_\mathcal{T}(\rho_i,Z_i^1)=l(i)d_{W^i}(0,U_i)$. Hence we
can also write $D_i=d_{W^i}(0,U_i)$, and so, for each $n$,
$(D_i)_{i\in\Sigma_n}$ is a collection of independent random
variables, independent of $\mathcal{F}_n$. Moreover, the random
variables $(D_i)_{i\in\Sigma_*}$ are identically distributed as
$D_\emptyset$, which represents the height of a $\mu$-random vertex in
$\mathcal{T}$. It is known that such a random variable has mean
$\sqrt{\pi/8}$, and finite variance (see \cite{Aldous2}, Section
3.3). Finally, we have the following recursive relationship
\begin{equation}D_i=w(i1)D_{i1}+w(i2)D_{i2},\label{drecurs}
\end{equation}
which may be deduced by decomposing the path from $\rho_i$ to $Z_i^1$
at $b(\rho_i,Z_i^1,Z_i^2)$, and applying parts (c) and (d) of
Proposition \ref{decompprop}.

\section{Self-similar dendrite in $\mathbb{R}^2$}\label{selfsimsec}

\begin{figure}[ht]
\centering
\scalebox{0.7}{
\includegraphics{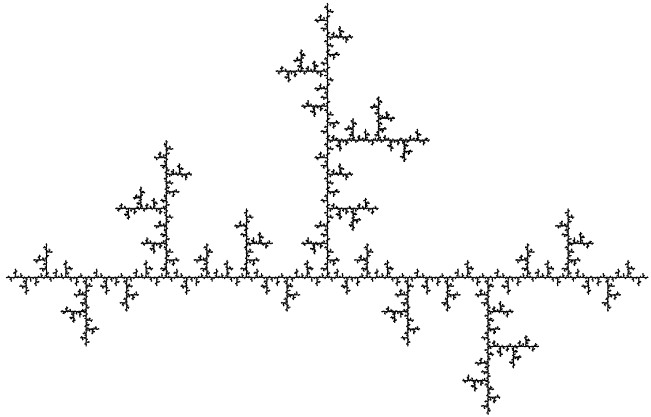}}
\caption{Self-similar dendrite.}
\label{ssd}
\end{figure}

The subset of $\mathbb{R}^2$ to which we will map the continuum random tree is a simple self-similar fractal, and is described as the fixed point of a collection of contraction maps. In particular, for $(x,y)\in\mathbb{R}^2$, set
\[F_1(x,y):=\frac{1}{2}(1-x,y),\hspace{10pt}F_2(x,y):=\frac{1}{2}(1+x,-y),\]
\[F_3(x,y):=\left(\frac{1}{2}+cy,cx\right),\]
where $c\in(0,1/2)$ is a constant, and define $T$ to be the unique
non-empty compact set satisfying $A=\bigcup_{i=1}^3F_i(A)$. The
existence and uniqueness of $T$, which is shown in Figure \ref{ssd},
is guaranteed by an extension of the usual contraction principle for
metric spaces, see \cite{Kigami}, Theorem 1.1.4. For a wide class of
self-similar fractals, which includes $T$, there is now a
well-established approximation procedure for defining an intrinsic
Dirichlet form and associated resistance metric on the relevant space,
see \cite{Barlow} and \cite{Kigami} for details. However, to capture
the randomness of the continuum random tree, we will need to randomise
this construction, and it is to describing how this is done that this
section is devoted.

The scaling factors that will be useful in defining a sequence of
compatible Dirichlet forms on subsets of $T$ will be the family
$(w(i))_{i\in\Sigma_*\backslash\{\emptyset\}}$, as defined in the
previous section. Although we would like to simply replace the
deterministic scaling factors that are used in the method of
\cite{Kigami} with this collection of random variables, following this
course of action would result in a sequence of non-compatible
quadratic forms, and taking limits would not be straightforward. To
deal with the offending tail fluctuations caused by using random
scaling factors, we introduce another collection of random variables
\begin{equation}\label{ridef}
R_i:=\lim_{n\rightarrow\infty}\sum_{j\in\{1,2\}^n}\frac{l(ij)}{l(i)},\hspace{20pt}i\in\Sigma_*,
\end{equation}
which we shall term resistance perturbations. Clearly these are
identically distributed, and, by appealing to the independence
properties of $(w(i))_{i\in\Sigma_*\backslash\{\emptyset\}}$, various
questions regarding the convergence and distribution of the
$(R_i)_{i\in\Sigma_*}$ may be answered by standard multiplicative
cascade techniques. Consequently we provide only a brief explanation
and suitable references for the proof of the following
result. Crucially, part (d) reveals an important identity between the
resistance perturbations and the family $(D_i)_{i\in\Sigma_*}$, which
was defined from the continuum random tree.

{\lem \label{resprop} (a) $\mathbf{P}$-a.s., the limit at
(\ref{ridef}) exists in $(0,\infty)$ for every $i\in\Sigma_*$.\\ (b)
$\mathbf{E}R_\emptyset =1$, and $\mathbf{E}R_\emptyset^d<\infty$ for
every $d\geq 0$.\\ (c) $\mathbf{P}$-a.s., for every $i\in\Sigma_*$,
the identity $R_i=w(i1)R_{i1}+w(i2)R_{i2}$ holds.\\ (d)
$\mathbf{P}$-a.s., $(R_i)_{i\in\Sigma_*}\equiv (HD_i)_{i\in\Sigma_*}$,
where $H:=\sqrt{8/\pi}$.}
\begin{proof}
The finite limit result of (a) and part (b) are immediate applications
of Theorem 2.0 of \cite{Liu}. Part (c) is immediate from the
definition of $(R_i)_{i\in\Sigma_*}$. Using the identical distribution
of the family of resistance perturbations, part (c) implies that
$\mathbf{P}(R_i=0)=\mathbf{P}(R_i=0)^2$. Since $\mathbf{E}R_i=1$, it
follows that $\mathbf{P}(R_i=0)=0$, which completes the proof of
(a). Checking the $\mathbf{P}$-a.s. equivalence of (d) is straightforward. First, from
an elementary application of a conditional version of Chebyshev's
inequality it may be deduced that, for each $i$,
\begin{eqnarray*}
\lefteqn{\mathbf{P}\left(\left|HD_i-\sum_{j\in\{1,2\}^n}\frac{l(ij)}{l(i)}\right|>\lambda\vline\:\mathcal{F}_{n+|i|}\right)}\\
&\leq&\lambda^{-2}\mathrm{Var}\left(\sum_{j\in\{1,2\}^n}\frac{l(ij)}{l(i)}(HD_{ij}-1)\vline\:\mathcal{F}_{n+|i|}\right)\\
&\leq&H^{2}\lambda^{-2}\sum_{j\in\{1,2\}^n}\frac{l(ij)^2}{l(i)^2}\mathrm{Var}D_{\emptyset}
\end{eqnarray*}
where we have also used the facts that $\mathbf{E}(HD_i)=1$, the identity at
(\ref{drecurs}) and the independence properties of the relevant random variables. Taking expectations yields
\[\mathbf{P}\left(\left|HD_i-\sum_{j\in\{1,2\}^n}\frac{l(ij)}{l(i)}\right|>\lambda\right)\leq
H^2\lambda^{-2}(\mathbf{E}(\Delta_1+\Delta_2))^n\mathrm{Var}D_\emptyset.\]
As remarked in the previous section, $D_\emptyset$ has finite
variance. Furthermore, a simple symmetry argument yields that the
expectation in the right hand side is precisely 2/3. Hence the sum of
probabilities over $n$ is finite, and applying a Borel-Cantelli
argument yields the result.\end{proof}

The sequence of vertices upon which we will define our Dirichlet forms
will be that which is commonly used for a p.c.f.s.s. fractal, see
\cite{Kigami} for more examples. Thus we shall not detail the reason
for the choice, but start by simply stating that the boundary of $T$
may be taken to be the two point set $V^0:=\{(0,0), (1,0)\}$. Our
initial Dirichlet form is defined by
\[D(f,f):=\sum_{x,y\in V^0,\:x\neq
y}H(f(x)-f(y))^2,\hspace{20pt}\forall f\in C(V^0),\]
where, for a countable set, $A$, we denote $C(A):=\{f:A\rightarrow
\mathbb{R}\}$. The constant $H$ is defined as in the previous lemma
and is necessary to achieve the correct scaling in the metric we shall
later define. We now introduce an increasing family of subsets of $T$
by setting $V^n:=\bigcup_{i\in\Sigma_n}F_i(V^0)$, where for
$i\in\Sigma_*$, $F_i:=F_{i_1}\circ\dots\circ F_{i_{|i|}}$. By defining
\[\mathcal{E}^n(f,f):=\sum_{i\in\Sigma_n}\frac{1}{l(i)R_i}D(f\circ
F_i,f\circ F_i),\hspace{20pt}\forall f\in C(V^n),\]
we obtain Dirichlet forms on each of the appropriate finite subsets of $T$,
$\mathbf{P}$-a.s. By applying the identity of Lemma \ref{resprop}(c),
it is straightforward to check that the family $(V^n,\mathcal{E}^n)$
is compatible in the sense that the trace of $\mathcal{E}^{n+1}$ on $V^n$ is precisely $\mathcal{E}^n$ for each $n$ (cf. \cite{Kigami}, Definition 2.2.1), and
from this fact we may take a limit in a sensible way. Specifically,
let
\[\mathcal{E}'(f,f):=\lim_{n\rightarrow\infty}\mathcal{E}^n(f,f),\hspace{20pt}\forall
f\in\mathcal{F}',\]
where $\mathcal{F}'$ is the set of functions on
the countable set $V^*:=\bigcup_{n\geq 0}V^n$ for which this limit
exists finitely. Note that we have abused notation slightly by using
the convention that if a form $\mathcal{E}$ is defined for functions
on a set $A$ and $f$ is a function defined on $B\supseteq A$, then we
write $\mathcal{E}(f,f)$ to mean $\mathcal{E}(f|_A, f|_A)$.

The quadratic form $(\mathcal{E}', \mathcal{F}')$ is actually a
resistance form (see \cite{Kigami}, Definition 2.3.1), and we can use
it to define a (resistance) metric $R'$ on $V^*$ using a formula
analogous to (\ref{dtrecover})
\[ R'(x,y)^{-1} = \inf\{\mathcal{E}'(f,f):f\in\mathcal{F}',
f(x)=0,f(y)=1\}, \]
for $x,y\in V^*$, $x\neq y$, and setting $R'(x,x)=0$.
We note that for sets of the form $F_i(V_0)$ with $i\in\Sigma_*$ we
have
\begin{equation}\label{edgeres}
R'(F_i(0,0),F_i(1,0))=\frac{l(i)R_i}{H}.
\end{equation}

To prove that this metric may be extended to
$T$ in a natural way (at least $\mathbf{P}$-a.s.) requires a similar
argument to the deterministic case, and so we omit the full details
here. The most crucial fact that is needed is the following:
\begin{equation}\label{diamdecay}
\lim_{n\rightarrow\infty}\sup_{i\in\Sigma_n}\mathrm{diam}_{R'}F_i(V^*)=0,
\hspace{20pt}\mathbf{P}\mbox{-a.s.},
\end{equation}
where, in general, $\mathrm{diam}_{d}(A)$ represents the diameter of a
set $A$ with respect to a metric $d$.
The proof follows the chaining argument of \cite{Barlow},
Proposition 7.10, and full details of the proof of the following
Proposition can be found in \cite{thesis}.
{\propn \label{rextend} There exists a unique metric $R$ on $T$ such
  that $(T,R)$ is the completion of $(V^*, R')$,
  $\mathbf{P}$-a.s. Moreover, the topology induced upon $T$ by $R$ is
  the same as that induced by the Euclidean metric, $\mathbf{P}$-a.s.}

To complete this section, we introduce the natural stochastic
self-similar measure on $T$, and note that
$(\mathcal{E}',\mathcal{F}')$ may be extended to a Dirichlet form on
the corresponding $L^2$ space. In particular, by proceeding exactly as
in the deterministic case, see \cite{Kigami}, Section 1.4, it is
possible to prove that, $\mathbf{P}$-a.s., there exists a unique
non-atomic Borel probability measure, $\mu^T$ say, on $(T,R)$ that
satisfies
\begin{equation}\label{selfsimmeas}
\mu^T(F_i(T))=l(i)^2,\hspace{20pt}\forall i\in\Sigma_*.
\end{equation}
Again, full details of this result are given in \cite{thesis}. If we
extend $(\mathcal{E}',\mathcal{F}')$ in the natural way by setting
$\mathcal{E}(f,f):=\mathcal{E}'(f,f)$, for $f \in \mathcal{F}:=\{f\in
C(T):\:f|_{V^*}\in\mathcal{F}'\}$, where we use $C(T)$ to represent
the continuous functions on $T$ (with respect to the Euclidean metric
or $R$), then the following result holds (for a proof, see
\cite{thesis}).

{\propn \label{rdescrp} $\mathbf{P}$-a.s., $(\mathcal{E},\mathcal{F})$
  is a local, regular Dirichlet form on $L^2(T,\mu^T)$ and, moreover,
  it may be associated with the metric $R$ through
\[ R(x,y)^{-1} = \inf\{\mathcal{E}(f,f): f\in\mathcal{F},
  f(x)=0,f(y)=1\}. \]}

\section{Equivalence of measure-metric spaces}\label{isosec}

In this section, we demonstrate how the decomposition of the continuum
random tree presented in Section \ref{decompsec} allows us to define
an isometry from the continuum random tree to the random self-similar
dendrite, $(T,R)$, described in the previous section. An important
consequence of the decomposition is that it allows us to label points
in $\mathcal{T}$ using the shift space of infinite sequences,
$\Sigma:=\{1,2,3\}^\mathbb{N}$. The following lemma defines the
projection $\pi_{\mathcal{T}}:\Sigma\rightarrow\mathcal{T}$ that we
will use, which is analogous to the well-known projection map for
self-similar fractals, see \cite{Barlow}, Lemma 5.10. We include the
result for the corresponding projection $\pi_T:\Sigma\rightarrow T$ to
allow us to introduce the necessary notation, and provide a direct
comparison of the two maps. Henceforth, we shall use the notation
$T_i:=F_i(T)$, for $i\in\Sigma_*$, and assume that $\Sigma$ is endowed with the usual ultra-metric topology generated by the sets $\{ij:j\in\Sigma\}$, $i\in\Sigma_*$.

{\lem \label{5.1} (a) There exists a map $\pi_{T}:\Sigma\rightarrow
{T}$ such that $\pi_{{T}}\circ\sigma_i(\Sigma)={T}_i$, for every
$i\in\Sigma_*$, where $\sigma_i:\Sigma\rightarrow\Sigma$ is defined by
$\sigma_i(j)=ij$ for $j\in\Sigma$. Furthermore, this map is
continuous, surjective and unique.\\ (b) $\mathbf{P}$-a.s., there
exists a map $\pi_\mathcal{T}:\Sigma\rightarrow \mathcal{T}$ such that
$\pi_{\mathcal{T}}\circ\sigma_i(\Sigma)=\mathcal{T}_i$, for every
$i\in\Sigma_*$, where $\sigma_i$ is defined as in (a). Furthermore,
this map is continuous, surjective and unique.}
\begin{proof} Part (a) is proved in \cite{Barlow} and \cite{Kigami},
  so we prove only (b). $\mathbf{P}$-a.s., for each $i\in\Sigma$, the
  sets in the collection $(\mathcal{T}_{i|n})_{n\geq 0}$ are compact,
  non-empty subsets of $(\mathcal{T},d_\mathcal{T})$, and by
  Proposition \ref{decompprop}(a), the sequence is decreasing. Hence,
  to show that $\cap_{n\geq0}\mathcal{T}_{i|n}$ contains exactly one
  point for each $i\in\Sigma$, $\mathbf{P}$-a.s., it will suffice to
  demonstrate that, $\mathbf{P}$-a.s.,
\begin{equation}\label{diamshrink}
\lim_{n\rightarrow\infty}\sup_{i\in\Sigma_n}\mathrm{diam}_{d_\mathcal{T}}\mathcal{T}_i=
0.
\end{equation}
From Proposition \ref{decompprop}(b), we have that
$\mathrm{diam}_{d_\mathcal{T}}\mathcal{T}_i=l(i)\mathrm{diam}_{d_{\tilde{\mathcal{T}}_i}}\tilde{\mathcal{T}}_i$. Using
the similarity that this implies, the above result may be proved in
the same way as (\ref{diamdecay}). To enable us to apply this
argument, we note that
$\mathrm{diam}_{d_{\tilde{\mathcal{T}}_i}}\tilde{\mathcal{T}}_i\leq2\sup_{t\in[0,1]}W^i_t$. The
upper bound here is simply twice the maximum of a normalised Brownian
excursion, and has finite positive moments of all orders as required
(see \cite{Aldous2}, for example).

Using the result of the previous paragraph, it is
$\mathbf{P}$-a.s. possible to define a map
$\pi_\mathcal{T}:\Sigma\rightarrow\mathcal{T}$ such that, for
$i\in\Sigma$,
$\{\pi_\mathcal{T}(i)\}=\bigcap_{n\geq0}\mathcal{T}_{i|n}$. That
$\pi_\mathcal{T}$ satisfies the claims of the lemma, and is the unique
map to do so, may be proved in exactly the same way as in the
self-similar fractal case.
\end{proof}

Heuristically, the isometry that we will define between the two
dendrites under consideration can be thought of as simply
``$\varphi=\pi_T\circ\pi_\mathcal{T}^{-1}$''. However, to introduce
the map rigorously, so that it is well-defined, we first need to prove
some simple, but fundamental, results about the geometry of the sets
and the maps $\pi_T$ and $\pi_\mathcal{T}$. From here on we use the notation $\dot{2}=222\dots$.

{\lem \label{lemmadrop} $\mathbf{P}$-a.s.,\\
(a) $\pi_\mathcal{T}^{-1}(\rho_{k1})=\{k11\dot{2},k21\dot{2},k31\dot{2}\}$, for all $k\in\Sigma_*$.\\
(b) For every $i,j\in\Sigma$, $\pi_\mathcal{T}(i)=\pi_{\mathcal{T}}(j)$ if and only if $\pi_T(i)=\pi_T(j)$.}
\begin{proof} The proof we give holds on the $\mathbf{P}$-a.s. set for
  which the decomposition of $\mathcal{T}$ and the definition of
  $\pi_\mathcal{T}$ is possible. Recall that
  $\rho_{k1}=b(\rho_k,Z_k^1,Z_k^2)$. For this branch-point to equal
  $\rho_{k}$ or $Z_k^1$, we would require at least two of its
  arguments to be equal, which happens with zero probability. Thus
  $\rho_{k1}\in\mathcal{T}_{k}\backslash\{\rho_k,Z_k^1\}$, and so
  Proposition \ref{decompprop}(f) implies that if
  $\pi_\mathcal{T}(i)=\rho_{k1}$ for some $i\in\Sigma$, then
  $i||k|=k$. Given this fact, it is elementary to apply the defining
  property of $\pi_\mathcal{T}$ and the results about $\rho_i$ and
  $Z_i^1$ that were deduced in Proposition \ref{decompprop} to deduce
  that part (a) of this lemma also holds. It now remains to prove part
  (b).

Fix $i,j\in\Sigma$, $i\neq j$, and let $m$ be the unique integer
satisfying $i|m=j|m$ and $i_{m+1}\neq j_{m+1}$. Furthermore, define
$k=i_1\dots i_m\in\Sigma_*$. Now by standard arguments for
p.c.f.s.s. fractals (see \cite{Kigami}, Proposition 1.2.5 and the
subsequent remark) we have that $\pi_T(i)=\pi_T(j)$ implies that
$\sigma^m(i),\sigma^m(j)\in\mathcal{C}$, where $\mathcal{C}$ is the
critical set for the self-similar structure, $T$, as defined in
\cite{Kigami}, Definition 1.3.4. Here, we use the notation $\sigma$ to
represent the shift map, which is defined by
$\sigma(i)=i_2i_3\dots$. Note that it is elementary to calculate that
$\mathcal{C}=\{11\dot{2},21\dot{2},31\dot{2}\}$ for this
structure. Thus $i,j\in\{k11\dot{2},k21\dot{2},k31\dot{2}\}$, and so,
by part (a), $\pi_\mathcal{T}(i)=\rho_{k1}=\pi_\mathcal{T}(j)$, which
completes one implication of the desired result.

Now suppose $\pi_\mathcal{T}(i)=\pi_\mathcal{T}(j)$. From the
definition of $\pi_\mathcal{T}$, we have that
$\pi_\mathcal{T}(i)\in\mathcal{T}_{ki_{m+1}}$ and also
$\pi_\mathcal{T}(j)\in\mathcal{T}_{kj_{m+1}}$. Hence
$\pi_\mathcal{T}(i),\pi_\mathcal{T}(j)\in\mathcal{T}_{ki_{m+1}}\cap\mathcal{T}_{kj_{m+1}}=\{\rho_{k1}\}$,
where we use (\ref{intersection}) to deduce the above equality. In
particular, this allows us to apply part (a) to deduce that
$i,j\in\{k11\dot{2},k21\dot{2},k31\dot{2}\}$. Applying the shift map
to this $m$ times yields $\sigma^m(i),\sigma^m(j)\in\mathcal{C}$. It
is easy to check that $\pi_T(\mathcal{C})$ contains only the single
point $(\frac{1}{2},0)$. Thus $\pi_T(i)=F_k\circ
\pi_T(\sigma^m(i))=F_k\circ\pi_T(\sigma^m(j))=\pi_T(j)$, which
completes the proof.
\end{proof}

We are now able to define the map $\varphi$ precisely on a $\mathbf{P}$-a.s. set by\begin{eqnarray*}
\varphi:\mathcal{T}&\rightarrow& T\\
x&\mapsto & \pi_T(i),\hspace{20pt}\mbox{for any $i\in \Sigma$ with $\pi_\mathcal{T}(i)=x$.}
\end{eqnarray*}
By part (b) of the previous lemma, this is a well-defined
injection. Furthermore, since $\pi_T$ is surjective, so is
$\varphi$. Hence we have constructed a bijection from $\mathcal{T}$ to
$T$ and it remains to show that it is also an isometry. We start by
checking that $\varphi$ is continuous, which will enable us to deduce
that it maps geodesic paths in $\mathcal{T}$ to geodesic paths in
$T$. However, before we proceed with the lemma, we introduce the
following notation for $x\in\mathcal{T}$, $n\geq 0$,
\[\mathcal{T}_n(x):=\bigcup\{\mathcal{T}_i:\:i\in\Sigma_n,\:x\in\mathcal{T}_i\}.\]
Define $(T_n(x))_{x\in T,n\geq 0}$ similarly, replacing
$\mathcal{T}_i$ with $T_i$ in the above definition where
appropriate. From the properties $\pi_T(i\Sigma)=T_i$,
$\pi_\mathcal{T}(i\Sigma)=\mathcal{T}_i$, and the definition of
$\varphi$, it is straightforward to deduce that
\begin{equation}\label{setmap}
\varphi(\mathcal{T}_i)=T_i,\hspace{20pt}\forall i\in\Sigma_*,
\end{equation}
on the $\mathbf{P}$-a.s. set that we can define all the relevant objects.

{\lem $\mathbf{P}$-a.s., $\varphi$ is a continuous map from $(\mathcal{T},d_\mathcal{T})$ to $(T,R)$.}
\begin{proof} By \cite{Kigami}, Proposition 1.3.6, for each $x\in T$,
  the collection $(T_n(x))_{n\geq 0}$ is a base of neighbourhoods of
  $x$ with respect to the Euclidean metric on $\mathbb{R}^2$. Since,
  by Proposition \ref{rextend}, $R$ is topologically equivalent to
  this metric, $\mathbf{P}$-a.s., then the same is true when we
  consider the collections of neighbourhoods with respect to $R$,
  $\mathbf{P}$-a.s. Similarly, we may use (\ref{diamshrink}),
  $\mathbf{P}$-a.s., to imitate the proofs of these results to deduce
  that $\mathbf{P}$-a.s., for each $x\in\mathcal{T}$, the collection
  $(\mathcal{T}_n(x))_{n\geq 0}$ is a base of neighbourhoods of $x$
  with respect to $d_\mathcal{T}$.

The remaining argument applies $\mathbf{P}$-a.s. Let $U$ be an open
subset of $(T,R)$ and $x\in\varphi^{-1}(U)$. Define $y=\varphi(x)\in
U$. Now, since $U$ is open, there exists an $n$ such that
$T_n(y)\subseteq U$. Also, by (\ref{setmap}), for each $i\in\Sigma_n$,
we have that $x\in \mathcal{T}_i$ implies that $y\in T_i$. Hence
\[\varphi(\mathcal{T}_n(x))=\varphi\left(\cup_{i\in\Sigma_n,\:x\in\mathcal{T}_i} \mathcal{T}_i\right)\subseteq\cup_{i\in\Sigma_n,\:y\in T_i} T_i=T_n(y)\subseteq U.\]
Consequently, $\mathcal{T}_n(x)\subseteq\varphi^{-1}(U)$. Since
$\mathcal{T}_n(x)$ is a $d_\mathcal{T}$-neighbourhood of $x$ it
follows that $\varphi^{-1}(U)$ is open in
$(\mathcal{T},d_\mathcal{T})$. The lemma follows.
\end{proof}

We are now ready to proceed with the main result of this section. In
the proof, we will use the notation
$\gamma_{xy}^{\mathcal{T}}:[0,1]\rightarrow\mathcal{T}$ to denote a
geodesic path (continuous injection) from $x$ to $y$, where $x$ and
$y$ are points in the dendrite $\mathcal{T}$. Clearly, because
$\varphi$ is a continuous injection,
$\varphi\circ\gamma_{xy}^{\mathcal{T}}$ describes a geodesic path from
$\varphi(x)$ to $\varphi(y)$ in $T$.

{\thm $\mathbf{P}$-a.s., the map $\varphi$ is an isometry, and the metric spaces $(\mathcal{T},d_\mathcal{T})$ and $(T,R)$ are isometric.}
\begin{proof} Obviously, the second statement of the theorem is an immediate consequence of the first. The following argument, in which we demonstrate that $\varphi$ is indeed an isometry, holds $\mathbf{P}$-a.s. Given $\varepsilon>0$, by (\ref{diamdecay}) and (\ref{diamshrink}), we can choose an $n\geq 1$ such that
\[\sup_{i\in\Sigma_n}\mathrm{diam}_{d_\mathcal{T}}\mathcal{T}_i,\:\sup_{i\in\Sigma_n}\mathrm{diam}_R
T_i<\frac{\varepsilon}{4}.\] Now, fix $x,y\in\mathcal{T}$, define
$t_0:=0$ and set
\[t_{m+1}:=\inf\{t>t_m\: :\:\gamma_{xy}^\mathcal{T}(t)\not\in\mathcal{T}_n(\gamma^\mathcal{T}_{xy}(t_m))\},\]
where $\inf\emptyset:=1$. We will also denote $x_m:=\gamma_{xy}^\mathcal{T}(t_m)$. Since, for each $x'\in\mathcal{T}$, the collection $(\mathcal{T}_n(x'))_{n\geq 0}$ forms a base of neighbourhoods of $x'$, we must have that $t_{m-1}<t_m$ whenever $t_{m-1}<1$. We now claim that for any $m$ with $t_{m-1}<1$ there exists a unique $i(m)\in\Sigma_n$ such that
\begin{equation}\label{claimjop}
\gamma_{xy}^\mathcal{T}(t)\in\mathcal{T}_{i(m)},\hspace{20pt}t_{m-1}\leq t\leq t_m.
\end{equation}
Let $m$ be such that $t_{m-1}<1$. By the continuity of $\gamma_{xy}^\mathcal{T}$, we have that $x_m\in\mathcal{T}_n(x_{m-1})$, and hence there exists an $i(m) \in\Sigma_n$ such that $x_{m-1},x_m\in\mathcal{T}_{i(m)}$. Clearly, the image of $\gamma^\mathcal{T}_{xy}$ restricted to $t\in[t_{m-1},t_m]$ is the same as the image of $\gamma_{x_{m-1}x_{m}}^\mathcal{T}$, which describes the unique path in $\mathcal{T}$ from $x_{m-1}$ to $x_m$. Note also that $\mathcal{T}_{i(m)}$ is a path-connected subset of $\mathcal{T}$, and so the path from $x_{m-1}$ to $x_m$ lies in $\mathcal{T}_{i(m)}$. Consequently, the set $\gamma_{xy}([t_{m-1},t_m])$ is contained in $\mathcal{T}_{i(m)}$. Thus to prove the claim at (\ref{claimjop}), it remains to show that $i(m)$ is unique. Suppose that there exists $j\in\Sigma_n$, $j\neq i(m)$ for which the inclusion at (\ref{claimjop}) holds. Then the uncountable set $\gamma_{xy}^\mathcal{T}([t_{m-1},t_m])$ is contained in $\mathcal{T}_{i(m)}\cap\mathcal{T}_j$, which, by Proposition \ref{decompprop}(f), contains at most two points. Hence no such $j$ can exist.

Now assume that $m_1<m_2$ and that $t_{m_2-1}<1$. Suppose that $i({m_1})=i({m_2})$, then $x_{m_1-1},x_{m_2}\in\mathcal{T}_{i({m_1})}$ By a similar argument to the previous paragraph, it follows that $\gamma_{xy}^\mathcal{T}([t_{m_1-1},t_{m_2}])\subseteq\mathcal{T}_{i({m_1})}$. By definition, this implies that $t_{m_1}\geq t_{m_2}$, which cannot be true. Consequently, we must have that $i({m_1})\neq i({m_2})$. Since $\Sigma_n$ is a finite set, it follows from this observation that $N:=\inf\{m:\:t_{m}=1\}$ is finite, and moreover, the elements of $(i(m))_{m=1}^N$ are distinct.

The conclusion of the previous paragraph provides us with a useful decomposition of the path from $x$ to $y$, which we will be able to use to complete the proof. The fact that $d_\mathcal{T}$ is a shortest path metric allows us to write $d_\mathcal{T}(x,y)=\sum_{m=1}^Nd_{\mathcal{T}}(x_{m-1},x_m)$. For $m\in\{2,\dots,N-1\}$, we have that $i({m})\neq i(m+1)$, and so by applying Proposition \ref{decompprop}(f), we can deduce that $x_{m}\in\mathcal{T}_{i(m)}\cap\mathcal{T}_{i({m+1})}\subseteq\{\rho_{i(m)},Z_{i(m)}^1\}$. Similarly, we have $x_{m-1}\in\mathcal{T}_{i(m-1)}\cap\mathcal{T}_{i(m)}\subseteq\{\rho_{i(m)},Z_{i(m)}^1\}$. Thus, by the injectivity of $\gamma_{xy}^\mathcal{T}$, we must have that $\{x_{m-1},x_{m}\}=\{\rho_{i(m)},Z_{i(m)}^1\}$, which implies $d_\mathcal{T}(x_{m-1},x_m)=d_{\mathcal{T}}(\rho_{i(m)},Z_{i(m)})=l(i(m))D_{i(m)}$. Hence we can conclude that
\begin{equation}\label{dtsplit}
d_\mathcal{T}(x,y)-\sum_{m=2}^{N-1}l(i(m))D_{i(m)}=d_{\mathcal{T}}(x_0,x_1)+d_\mathcal{T}(x_{N-1},x_N).
\end{equation}

As remarked before this lemma, $\varphi\circ\gamma_{xy}^\mathcal{T}$ is a geodesic path from $\varphi(x)$ to $\varphi(y)$. Thus the shortest path property of $R$ allows us to write
\begin{equation}\label{rsplit}
R(\varphi(x),\varphi(y))=\sum_{m=1}^N R(\varphi(x_{m-1}),\varphi(x_m)).
\end{equation}
Let $m\in\{2,\dots,N-1\}$. By applying $\varphi$ to the expression for $\{x_{m-1},x_m\}$ that was deduced above, we obtain that $\{\varphi(x_{m-1}),\varphi({x_m})\}=\{\varphi(\rho_{i(m)}),\varphi(Z_{i(m)}^1)\}$. Now, part (a) of Lemma \ref{lemmadrop} implies that
\[\varphi(\rho_{i(m)})=\pi_T(k11\dot{2})=F_{k}(\pi_T(11\dot{2}))=F_k((\frac{1}{2},0))=F_{i(m)}((0,0)),\]
where $k:=i(m)|(|i(m)|-1)$. In Proposition \ref{decompprop}(d) it was shown that $Z_i^1=Z_{i2}^1$, for every $i\in\Sigma_*$. It follows that $i(m)\dot{2}\in\pi_{\mathcal{T}}^{-1}(Z_{i(m)}^1)$, and so
\[\varphi(Z^1_{i(m)})=\pi_T(i(m)\dot{2})=F_{i(m)}(\pi_T(\dot{2}))=F_{i(m)}((1,0)).\]
Thus $R(\varphi(x_{m-1}),\varphi(x_m))=R(F_{i(m)}((0,0)),F_{i(m)}((1,0)))$, and so from the expression at (\ref{edgeres}), we can deduce that $R(\varphi(x_{m-1}),\varphi(x_m))=\sqrt{{\pi}/{8}}l(i(m))R_{i(m)}$, which, by Lemma \ref{resprop}(d), is equal to $l(i(m))D_{i(m)}$. Substituting this into (\ref{rsplit}), and combining the resulting equation with the equality at (\ref{dtsplit}) yields
\begin{eqnarray*}
\lefteqn{|d_\mathcal{T}(x,y)-R(\varphi(x),\varphi(y))|\leq}\\
&&\sum_{m\in\{1,N\}}\left(d_\mathcal{T}(x_{m-1},x_m)+R(\varphi(x_{m-1}),\varphi(x_m))\right).
\end{eqnarray*}
Now, $x_0$ and $x_1$ are both contained in $\mathcal{T}_{i(1)}$, and so the choice of $n$ implies that $d_\mathcal{T}(x_0,x_1)<\varepsilon/4$. Furthermore, $\varphi(x_0)$ and $\varphi(x_1)$ are both contained in $\varphi(\mathcal{T}_{i(1)})=T_{i(1)}$, and so we also have $R(\varphi(x_{0}),\varphi(x_1))<\varepsilon/4$. Thus the summand with $m=1$ is bounded by $\varepsilon/2$. Similarly for $m=N$. Hence $|d_\mathcal{T}(x,y)-R(\varphi(x),\varphi(y))|<\varepsilon$. Since the choice of $x,y$ and $\varepsilon$ was arbitrary, the proof is complete.
\end{proof}

The final result that we present in this section completes the proof of the fact that $(\mathcal{T},d_{\mathcal{T}},\mu)$ and $(T,R,\mu^T)$ are equivalent measure-metric spaces, where we continue to use the notation $\mu^T$ to represent the stochastic self-similar measure on $(T,R)$, as defined in Section \ref{selfsimsec}.

{\thm $\mathbf{P}$-a.s., the probability measures $\mu$ and $\mu^T\circ\varphi$ agree on the Borel $\sigma$-algebra of $(\mathcal{T},d_\mathcal{T})$.}
\begin{proof} That both $\mu^T\circ\varphi$ and $\mu$ are non-atomic Borel probability measures on $(\mathcal{T},d_\mathcal{T})$, $\mathbf{P}$-a.s., is obvious. Recall from Proposition \ref{decompprop}(g) that $\mu(\mathcal{T}_i)=l(i)^2$, for every $i\in\Sigma_*$, $\mathbf{P}$-a.s. Furthermore, from the identities of (\ref{selfsimmeas}) and (\ref{setmap}), we also have $\mu^T\circ\varphi(\mathcal{T}_i)=\mu^T(T_i)=l(i)^2$, for every $i\in\Sigma_*$, $\mathbf{P}$-a.s. The result is readily deduced from these facts.
\end{proof}

\section{Spectral asymptotics}\label{specsec}

Due to the construction of the natural Dirichlet form on the continuum random tree from the natural metric on the space, the results of the previous section imply that the spectrum of $(\mathcal{E}_\mathcal{T},\mathcal{F}_\mathcal{T},\mu)$ is $\mathbf{P}$-a.s. identical to that of $(\mathcal{E},\mathcal{F},\mu^T)$, the random Dirichlet form and self-similar measure on $T$, as defined in Section \ref{selfsimsec}. Consequently, to deduce the results of the introduction, it will suffice to show that the analogous results hold for $(\mathcal{E},\mathcal{F},\mu^T)$, which is possible using techniques developed for related self-similar fractals. For this argument, it will be helpful to apply various decomposition and comparison inequalities for the Dirichlet and Neumann eigenvalues associated with this Dirichlet form, and we shall start by introducing these.

To define the Dirichlet eigenvalues for $(\mathcal{E},\mathcal{F},\mu^T)$, we first introduce the related Dirichlet form $(\mathcal{E}^D,\mathcal{F}^D)$ by setting
\[\mathcal{E}^D(f,f):=\mathcal{E}(f,f),\hspace{20pt}\forall f\in\mathcal{F}^D,\]
where
\[\mathcal{F}^D:=\{f\in\mathcal{F}:\:f|_{V^0}=0\}.\]
The Dirichlet eigenvalues of the original form, $(\mathcal{E},\mathcal{F},\mu^T)$, are then defined to be the eigenvalues of $(\mathcal{E}^D,\mathcal{F}^D, \mu^T)$. We shall use the title Neumann eigenvalues to refer to the usual eigenvalues of $(\mathcal{E},\mathcal{F},\mu^T)$, defined analogously to (\ref{evaluedef}).

Before continuing, note that the description of $R$ in Proposition \ref{rdescrp} easily leads to the well known inequality
\begin{equation}\label{resineq}
|f(x)-f(y)|^2\leq R(x,y)\mathcal{E}(f,f),\hspace{20pt}\forall x,y\in T,\:f\in\mathcal{F}.
\end{equation}
By applying this fact (and using $\|\cdot\|_p$ to represent the corresponding $L^p(T,\mu^T)$ norm), we find that, for $x\in T$, $f\in\mathcal{F}$,
\[|f(x)|^2 \leq 2\int_T (|f(x)-f(y)|^2+|f(y)|^2)d\mu\leq  2 \mathrm{diam}_R T \mathcal{E}(f,f)+2\|f\|^2_2,\]
and so, $\mathbf{P}$-a.s., $\|f\|_\infty^2\leq C(\mathcal{E}(f,f)+\|f\|^2_2)$, for some constant $C$. Combining this inequality with (\ref{resineq}), we can imitate the argument of \cite{KigLap}, Lemma 5.4, to deduce that the natural inclusion map from $(\mathcal{F},\mathcal{E}+\|\cdot\|_2^2)$ to $L^2(T,\mu^T)$ is a compact operator. It follows that the Dirichlet and Neumann spectra of $(\mathcal{E},\mathcal{F},\mu^T)$ are discrete, and so the associated eigenvalue counting functions, $N^D(\lambda)$ and $N^N(\lambda)$, are well-defined and finite for all $\lambda\in\mathbb{R}$.

From the definitions in the previous paragraph, we can easily see that
$N(\lambda)=N^N(\lambda)$, $\mathbf{P}$-a.s., and so, using the terminology
introduced above, the eigenvalues of $(\mathcal{E}_\mathcal{T},
\mathcal{F}_\mathcal{T}, \mu)$ may be thought of as Neumann eigenvalues. Of
course, this definition does not provide any justification for using the name
Neumann, so we will now give an explanation of why it is sensible to do so. Since we will not actually apply this interpretation, we only sketch the relevant results. Analogously to \cite{Kigami}, Definition 3.7.1, let $\mathcal{D}$ be the collection of functions $f\in C(T)$ such that there exists a function $g\in C(T)$ satisfying
\begin{equation}\label{limit}\lim_{n\rightarrow\infty}\max_{x\in V^n\backslash V^0}\left|\mu_{n,x}^{-1}\Delta_n f(x)-g(x)\right|=0,\end{equation}
where $\Delta_n$ is the discrete Laplacian on $V^n$ associated with $\mathcal{E}^n$, $\mu_{n,x}:=\int_T\psi_x^nd\mu$, and $\psi_x^n$ is the unique harmonic extension (with respect to $(\mathcal{E},\mathcal{F})$) of $\mathbf{1}_{\{x\}}$ from $V^n$ to $T$. For a function $f\in \mathcal{D}$ satisfying (\ref{limit}), we write $\Delta f =g$, so that $\Delta$ is essentially the limit operator of the rescaled discrete Laplacians $\Delta_n$. Furthermore, for $f\in\mathcal{D}$, we can also
define a function, $df$ say, with domain $V^0$, which represents the Neumann
derivative on the boundary of $T$ (similarly to \cite{Kigami}, Definition 3.7.6) by setting $(df)(x):=\lim_{n\rightarrow\infty}-\Delta_n u(x)$. By using a Green's function argument as in
the proof of \cite{Kigami}, Theorem 3.7.9, it is possible to deduce that the
Friedrichs extension of $\Delta$ on
$\mathcal{D}_D:=\{f\in\mathcal{D}:\:f|_{V^0}=0\}$ is precisely $\Delta_D$, the
Laplacian associated with $(\mathcal{E}^D, \mathcal{F}^D, \mu^T)$. Similarly,
the Friedrichs extension of $\Delta$ on
$\mathcal{D}_N:=\{f\in\mathcal{D}:\:(df)(x)=0,\:\forall x\in V^0\}$ is
$\Delta_N$, the Laplacian associated with $(\mathcal{E},\mathcal{F}, \mu^T)$.
Note that the construction of the relevant Green's function may be accomplished
more easily than in \cite{Kigami} by, instead of imitating the analytic
definition used there, applying a probabilistic definition, with $g(x,y)$ being
the Green's kernel for the Markov process associated with $(\mathcal{E},
\mathcal{F})$ killed on hitting $V^0$ (the existence of which follows from an
argument similar to that used in \cite{Kumagai}, Proposition 4.2).

Applying the relationships between the various operators introduced in the
previous paragraph (and also the continuity of the Green's function), we are
able to emulate the argument of \cite{Kigami}, Proposition 4.1.2, to deduce
that the eigenvalues of $(\mathcal{E}^D, \mathcal{F}^D, \mu^T)$ are precisely
the solutions to \[-\Delta u=\lambda u,\hspace{20pt}u|_{V^0}=0,\] for some
eigenfunction $u\in\mathcal{D}$. Furthermore, the eigenvalues of $(\mathcal{E},
\mathcal{F}, \mu^T)$ are precisely the solutions to \begin{equation}\label{Neu}
-\Delta u=\lambda u,\hspace{20pt}(du)|_{V^0}=0, \end{equation} for some
eigenfunction $u\in\mathcal{D}$. From these characterisations, it is clear that
the Dirichlet and Neumann eigenvalues of $(\mathcal{E}, \mathcal{F}, \mu^T)$
that we have defined are exactly the eigenvalues of $-\Delta$ with the usual
Dirichlet (zero function on boundary) and Neumann (zero derivative on boundary)
boundary conditions respectively, where the analytic boundary of $T$ is taken
to be $V^0$.

By mapping these results to the continuum random tree, we are able to deduce,
$\mathbf{P}$-a.s., the existence of a Laplace operator $\Delta_\mathcal{T}$ on
$\mathcal{T}$, and also a Neumann boundary derivative, so that the eigenvalues
of $(\mathcal{E},\mathcal{F},\mu)$ satisfy a result analogous to (\ref{Neu}).
In the continuum random tree setting, observe that the natural analytic
boundary is the two point set consisting of the root and one $\mu$-random
vertex, $\{\rho, Z^1_\emptyset\}$. Consequently, the results we prove also
demonstrate the Dirichlet spectrum corresponding to this boundary satisfies the
same asymptotics as the original (Neumann) spectrum. Another point of interest
is that by replicating the argument of \cite{Kigami}, Theorem 3.7.14, we are
able to uniquely solve the Dirichlet problem for Poisson's equation (with
respect to $\Delta_\mathcal{T}$) on the continuum random tree, again taking $\{\rho,Z^1_\emptyset\}$ as our boundary.

We now return to our main argument. From the construction of
$(\mathcal{E},\mathcal{F})$, it is possible to deduce the following
self-similar decomposition using the same proof as in Lemma 4.5 of
\cite{Hamasymp}.

{\lem $\mathbf{P}$-a.s., we have, for every $n\geq 1$,
\[\mathcal{E}(f,g)=\sum_{i\in\Sigma_n}\frac{1}{l(i)}\mathcal{E}_i(f\circ F_i,g\circ F_i),\hspace{20pt}\forall f,g\in\mathcal{F},\]
where $(\mathcal{E}_i)_{i\in\Sigma_n}$ are independent copies of $\mathcal{E}$, independent of $\mathcal{F}_n$.}
\bigskip

The operators of the above theorem each have a Dirichlet version, $\mathcal{E}_i^D$, defined in the same way as $\mathcal{E}^D$ was from $\mathcal{E}$. We shall denote by $N^D_i(\lambda)$ and $N^N_i(\lambda)$ the corresponding Dirichlet and Neumann eigenvalue counting functions.

{\lem \label{6.2} $\mathbf{P}$-a.s., we have, for every $\lambda>0$,
\begin{equation}\label{ineq}\sum_{i=1}^3N_i^D(\lambda w(i)^3)\leq N^D(\lambda)\leq N^N(\lambda)\leq \sum_{i=1}^3 N_{i}^N(\lambda w(i)^3),\end{equation}
and also $N^D(\lambda)\leq N^N(\lambda)\leq N^D(\lambda)+2$.}
\begin{proof} Since the proof of this result can be completed by repeating the argument of \cite{Hamasymp}, Lemma 5.1, we will only present a brief outline here. First, define a quadratic form $(\tilde{\mathcal{E}}^D,\tilde{\mathcal{F}}^D)$ by setting $\tilde{\mathcal{E}}^D=\mathcal{E}^D|_{\tilde{\mathcal{F}}^D\times\tilde{\mathcal{F}}^D}$, where $\tilde{\mathcal{F}}^D$ is the set $\{f\in\mathcal{F}^D:\:f|_{V^1}=0\}$. It is straightforward to check that $(\tilde{\mathcal{E}}^D,\tilde{\mathcal{F}}^D)$ is a local Dirichlet form on $L^2(T,\mu^T)$ and the natural inclusion map from $\tilde{\mathcal{F}}^D$ to $L^2(T,\mu^T)$ is compact, $\mathbf{P}$-a.s. Thus we can define the related eigenvalue counting function $\tilde{N}^D(\lambda)$ and, by \cite{Hamasymp}, Lemma 5.4, we have $\tilde{N}^D(\lambda)\leq N^D(\lambda)$ for all $\lambda$, $\mathbf{P}$-a.s. Now, fix $i\in\{1,2,3\}$ and suppose $f$ is an eigenfunction of $(\mathcal{E}_i^D,\mathcal{F}_i^D,\mu^T_i)$ with eigenvalue $\lambda w(i)^3$, where the domain $\mathcal{F}_i^D$ of $\mathcal{E}_i^D$ is defined analogously to $\mathcal{F}^D$ and $\mu^T_i:=\mu^T\circ F_i(\cdot) /\mu^T(T_i)$. If we set
\[g(x):=\left\{\begin{array}{ll}
          f\circ F_i^{-1} (x), & \mbox{for }x\in T_i, \\
          0 & \mbox{otherwise,}
        \end{array}\right.\]
then, by definition, we have that, for $h\in\tilde{\mathcal{F}}^D$,
\[\tilde{\mathcal{E}}^D(g,h)=\frac{1}{w(i)}\mathcal{E}_i^D(f,h\circ F_i)=\lambda w(i)^2\int_T f (h\circ F_i)d\mu^T_i=\lambda \int_T gh d\mu^T.\]
Thus $g$ is an eigenfunction of $(\tilde{\mathcal{E}}^D,\tilde{\mathcal{F}}^D,\mu^T)$ with eigenvalue $\lambda$. Hence it is clear that $\sum_{i=1}^3N_i^D(\lambda w(i)^3)\leq \tilde{N}^D(\lambda)$ for all $\lambda$, $\mathbf{P}$-a.s., which completes the proof of the left-hand inequality of (\ref{ineq}). The right-hand inequality of (\ref{ineq}) is proved using a similar decomposition of a suitable enlargement of $(\mathcal{E},\mathcal{F})$, see \cite{Hamasymp}, Proposition 5.2, for details. The remaining parts of the lemma are a simple application of Dirichlet-Neumann bracketing, see \cite{Hamasymp}, Lemma 5.4, for example.
\end{proof}

For the remainder of this section, we will continue to follow \cite{Hamasymp}, and proceed by defining a time-shifted general branching process, $X$. Although the results we shall prove will be in terms of $N^D$, the second set of inequalities in the above lemma imply that the asymptotics of $N^N$, and consequently $N$, are the same.

Define the functions $(\eta_i)_{i\in\Sigma_*}$ by, for $t\in\mathbb{R}$,
\[\eta_i(t):=N_i^D(e^t)-\sum_{j=1}^3N_{ij}^D(e^t w(ij)^3),\]
and let $\eta:=\eta_\emptyset$. Clearly, the paths of $\eta_i(t)$ are cadlag, and Lemma \ref{6.2} implies that the functions take values in $[0,6]$, $\mathbf{P}$-a.s. If we set $X_i(t):=N_i(e^t)$, and $X:=X_\emptyset$, then it is possible to check that the following evolution equation holds:
\begin{equation}\label{evo}
X(t)=\eta(t)+\sum_{i=1}^3X_i(t+3\ln w(i));
\end{equation}
and also that
\begin{equation}\label{alter}
X(t)=\sum_{i\in\Sigma_*}\eta_i(t+3\ln l(i)).
\end{equation}
The equation at (\ref{evo}) is particularly important, as it will allow us to use branching process and renewal techniques to obtain the results of interest.

We start by investigating the mean behaviour of $X$, and will now introduce the notation necessary to do this. Set $\gamma=2/3$, and define, for $t\in\mathbb{R}$,
\begin{equation}\label{mdef}
m(t):=e^{-\gamma t}\mathbf{E}X(t),\hspace{20pt}u(t):=e^{-\gamma t}\mathbf{E}\eta(t).
\end{equation}
Furthermore, define the measure $\nu$ by $\nu([0,t])=\sum_{i=1}^3 \mathbf{P}(w(i)^3\geq e^{-t})$, and let $\nu_\gamma$ be the measure that satisfies $\nu_\gamma(dt)=e^{-\gamma t}\nu(dt)$. Some properties of these objects are collected in the following lemma.

{\lem \label{props} (a) The function $m$ is bounded and  measurable, and $m(t)\rightarrow0$ as $t\rightarrow -\infty$.\\
(b) The function $u$ is in $L^1(\mathbb{R})$ and $u(t)\rightarrow0$ as $|t|\rightarrow \infty$.\\
(c) The measure $\nu_\gamma$ is a Borel probability measure on $[0,\infty)$, and the integral $\int_0^\infty t\nu_\gamma(dt)$ is finite.}
\begin{proof} A fact that may be deduced from (\ref{resineq}), and will be important in proving parts (a) and (b), is that $\mathbf{P}$-a.s., $\|f\|_2^2\leq \mathcal{E}(f,f)\mathrm{diam}_R T$, for every $f\in\mathcal{F}^D$. In particular, this implies that the bottom of the Dirichlet spectrum is bounded below by $(\mathrm{diam}_R T)^{-1}$, and consequently we must have $\eta(t)=0$ for $t<-\ln\mathrm{diam}_RT$, $\mathbf{P}$-a.s. Hence,
\begin{equation}\label{etazeronound}
\mathbf{E}\eta(t)\leq 6 \mathbf{P}(t\geq -\ln\mathrm{diam}_RT).
\end{equation}
Applying this result, the alternative representation of $X$ at (\ref{alter}), and the independence of $N_i^D$ and $\mathcal{F}_{|i|}$, we obtain
\begin{eqnarray*}
m(t)&=&\sum_{i\in\Sigma_*}e^{-\gamma t}\mathbf{E}\eta_i(t+3\ln l(i))\\
&\leq&6e^{-\gamma t}\mathbf{E}(\#\{i\in\Sigma_*:\:t+3\ln l(i)\geq-\ln\mathrm{diam}_{R'} T\}),
\end{eqnarray*}
where $R'$ is an independent copy of $R$. Applying standard branching process techniques to the process with particles $i\in\Sigma_*$, where $i\in\Sigma_*$ has offspring $ij$ at time $-\ln w(ij)$ after its birth, $j=1,2,3$, it is possible to show that $\mathbf{E}(\#\{i\in\Sigma_*:\:\ln l(i)\geq -t\})\leq Ce^{2t}$, for every $t\in\mathbb{R}$; the exponent 2 that arises is the Malthusian parameter for the relevant branching process. Thus, for $t\in\mathbb{R}$,
\[m(t)\leq 6C\mathbf{E}((\mathrm{diam}_RT)^\gamma).\]
Since $\mathrm{diam}_RT\buildrel{d}\over{=}\mathrm{diam}_{d_\mathcal{T}}\mathcal{T}$, and, as remarked in the proof of Lemma \ref{5.1}, $\mathrm{diam}_{d_\mathcal{T}}\mathcal{T}$ has finite positive moments, we are able to deduce that the right hand side of the above inequality is finite. Thus, $m$ is bounded. The measurability of $m$ follows from the fact that $X$ has cadlag paths, $\mathbf{P}$-a.s. To demonstrate the limit result, we recall the bound at (\ref{etazeronound}), which we apply to (\ref{alter}) to obtain
\[m(t)\leq \sum_{i\in\Sigma_*}6e^{-\gamma t}\mathbf{P}(l(i)^3\mathrm{diam}_{R'}T\geq e^{-t}),\]
where $R'$ is again an independent copy of $R$. Applying Markov's inequality to this expression, we find that, for $\theta>0$,
\begin{eqnarray*}
m(t)&\leq&\sum_{i\in\Sigma_*} 6 e^{(\theta-\gamma)t}\left(\mathbf{E}(w(i)^{3\theta})\right)^{|i|}\mathbf{E}((\mathrm{diam}_RT)^\theta)\\
&=& 6 e^{(\theta-\gamma)t} \mathbf{E}((\mathrm{diam}_RT)^\theta)\sum_{n\geq 0}3^n\left(\mathbf{E}(w(1)^{3\theta})\right)^{n}
\end{eqnarray*}
Taking $\theta>\gamma$, we have $\mathbf{E}(w(1)^{3\theta})<\mathbf{E}(w(1)^2)=\frac{1}{3}$, so the sum over $n$ is finite, as is the expectation involving $\mathrm{diam}_RT$. Consequently, the upper bound converges to zero as $t\rightarrow -\infty$, which completes the proof of (a).

That $u(t)$ is finite for $t\in\mathbb{R}$ follows from the fact that $\eta(t)$ is, and the measurability of $u$ is a result of $\eta$ having cadlag paths, $\mathbf{P}$-a.s. Observe that, for $t\geq 0$,
\begin{equation}\label{bound}
\mathbf{P}\left(\mathrm{diam}_RT>t\right)=\mathbf{P} \left(
\mathrm{diam}_{d_\mathcal{T}} \mathcal{T}>t\right)\leq
\mathbf{P}\left(\sup_{s\in[0,1]}W_s>\frac{t}{2}\right)\leq Ce^{-t^2/4},
\end{equation}
for some constant $C$, where the final inequality is obtained by
applying the exact distribution of the supremum of a normalised
Brownian excursion (see \cite{Aldous2}, Section 3.1). Thus, again applying (\ref{etazeronound}), we see that $u(t)$ is bounded above by $6 e^{-\gamma t}\left(1\wedge Ce^{-e^{-2t}/4}\right)$ for all $t$, which readily implies the remaining claims of (b).

Part (c) is easily deduced using simple properties of the Dirichlet
distribution of the triple $(w(1)^2,w(2)^2,w(3)^2)$.
\end{proof}

The importance of the previous lemma is that it allows us to apply the renewal
theorem to deduce the mean behaviour of $X$, with the precise result being
presented in the following proposition. Part (a) of Theorem \ref{second} is an
easy corollary of this.

{\propn \label{meanconv} The function $m$ converges as $t\rightarrow \infty$ to
the finite and non-zero constant \[m(\infty):=\frac{\int_{-\infty}^\infty
u(t)dt}{\int_0^\infty t\nu_\gamma(dt)}.\]}
\begin{proof} After multiplying by $e^{-\gamma t}$ and taking
  expectations, the equation at (\ref{evo}) may be rewritten, for
  $t\in\mathbb{R}$, \[m(t)=u(t)+\int_0^\infty m(t-s)\nu_\gamma(ds),\]
which is the double-sided renewal equation of \cite{Karlin}. The
  results that are proved about $m$, $u$ and $\nu_\gamma$ in Lemma
  \ref{props} mean that the conditions of the renewal theorem stated
  in \cite{Karlin} are satisfied, and the proposition follows from
  this. \end{proof}

To determine the $\mathbf{P}$-a.s. behaviour of $X$, and prove part (b) of
Theorem \ref{second}, the argument of \cite{Hamasymp}, Section 5, may be used.
Note that this method is in turn an adaptation of Nerman's results on the
almost-sure behaviour of general branching processes, see \cite{Nerman}. Since the steps
of our proof are almost identical to those of \cite{Hamasymp}, we shall omit
many of the details here. One point that should be highlighted, however, is
that in the proof of Lemma 5.7 of \cite{Hamasymp} there is an error, with
one of the relevant terms being omitted from consideration. We shall explain how to deal
with this term, and also correct the limiting procedure that should be used at
the end of the argument. For the purposes of the proof, we introduce the
following notation to represent a cut-set of $\Sigma_*$: for $t>0$,
\[\Lambda_t:=\{i\in\Sigma_*:\:-3\ln l(i)\geq t > -3\ln l(i|(|i|-1))\}.\] We
will also have cause to refer to the subset of $\Lambda_t$ defined by, for $t,c>0$,
\[\Lambda_{t,c}:=\{i\in\Sigma_*:\:-3\ln l(i)\geq t+c,\: t> -3\ln
l(i|(|i|-1))\}.\]

{\propn $\mathbf{P}$-a.s., we have
\[e^{-\gamma t}X(t)\rightarrow m(\infty),\hspace{20pt}\mbox{as }t\rightarrow 0,\]
where $m(\infty)$ is the constant defined in Proposition \ref{meanconv}.}
\begin{proof} First, we
truncate the characteristics $\eta_i$ by defining, for fixed $c>0$,
$\eta^c_i(t):=\eta_i(t)\mathbf{1}_{\{t<n_0c\}}$, where $n_0$ is an integer that
will be chosen later in the proof (we are using the term ``characteristic'' in the generalised sense of \cite{Nerman}, Section 7, to describe a random function that, if the characteristic is indexed by $i$, can depend on $(w(ij))_{j\in\Sigma_*}$). From these truncated characteristics, construct the
processes $X_i^c$, by
\[X_i^c(t):=\sum_{j\in\Sigma_*}\eta_{ij}^c(t+3\ln (l(ij)/l(i))),\]
and set $X^c:=X_\emptyset^c$. The corresponding discounted mean process is $m^c(t):=e^{-\gamma t}\mathbf{E}X^c(t)$, and this may be checked to converge to $m^c(\infty)\in(0,\infty)$ as
$t\rightarrow\infty$ using the argument of Proposition \ref{meanconv}. From a branching process decomposition of $X^c$, we can deduce the following bound for $n_1\geq n_0$, $n\in\mathbb{N}$,
\[|e^{-\gamma c(n+n_1)}X^c(c(n+n_1))-m^c(\infty)|\leq S_1(n,n_1)+S_2(n,n_1)+S_3(n,n_1),\]
where,
\begin{eqnarray*}
\lefteqn{S_1(n,n_1):=}\\
&&\left|\sum_{i\in\Lambda_{cn}\backslash\Lambda_{cn,cn_1}} \left(e^{-\gamma c (n+n_1)} X_i^c(c(n+n_1)+3\ln l(i))-l(i)^2 m^c(c(n+n_1)+3\ln l(i))\right)\right|,
\end{eqnarray*}
\[S_2(n,n_1):=\left|\sum_{i\in\Lambda_{cn}\backslash\Lambda_{cn,cn_1}}l(i)^2 m^c(c(n+n_1)+3\ln l(i))-m^c(\infty)\right|,\]
\[S_3(n,n_1):=e^{-\gamma c(n+n_1)}\sum_{i\in\Lambda_{cn,cn_1}}X_i^c(c(n+n_1)+3\ln l(i)).\]
The first two of these terms are dealt with in \cite{Hamasymp}, and using the arguments from that article, we have that, $\mathbf{P}$-a.s.,
\[\lim_{n_1\rightarrow\infty}\limsup_{n\rightarrow\infty} S_j(n,n_1)=0,\hspace{20pt}\mbox{for }j=1,2.\]

We now show how $S_3(n,n_1)$ decays in a similar fashion. First, introduce a set of characteristics, $\phi_i^{c,n_1}$, defined by
\[\phi_i^{c,n_1}(t):=\sum_{j=1}^3X_{ij}(0)\mathbf{1}_{\{-3\ln w(ij)>t+cn_1, \:t>0\}},\]and, for $t>0$, set
\[Y^{c,n_1}(t):=\sum_{i\in\Sigma_*}\phi_i^{c,n_1}(t+3\ln l(i)).\]
Note that from the definition of the cut-sets $\Lambda_{cn}$ and $\Lambda_{cn,cn_1}$ we can deduce that
\[Y^{c,n_1}(cn)=\sum_{i\in\Lambda_{cn,cn_1}}X_i(0)\geq e^{\gamma c (n+n_1)}S_3(n,n_1),\]
where for the second inequality we apply the monotonicity of the $X_i$s. Now, $Y^{c,n_1}$ is a branching process with random characteristics $\phi^{c,n_1}_i$, and we are able to check the conditions of the extension of \cite{Nerman}, Theorem 5.4, that is stated as \cite{Hamasymp}, Theorem 3.2, are satisfied. By applying this result, we find that $\mathbf{P}$-a.s.,
\[e^{-\gamma t}Y^{c,n_1}(t)\rightarrow \frac{\int_0^\infty e^{-\gamma t}\mathbf{E}\phi^{c,n_1}_\emptyset(t)dt}{\int_0^{\infty}t\nu_\gamma(dt)},\hspace{20pt}\mbox{as }t\rightarrow\infty.\]
It is obvious that $\mathbf{E}\phi^{c,n_1}_\emptyset(t)\leq 3\mathbf{E}X(0)\leq 3m(0)<\infty$, where $m$ is the function defined at (\ref{mdef}). Consequently, there exists a constant $C$ that is an upper bound for the above limit uniformly in $n_1$, and so $\mathbf{P}$-a.s.,
\[\lim_{n_1\rightarrow\infty}\limsup_{n\rightarrow\infty}S_3(n,n_1)\leq \lim_{n_1\rightarrow\infty} Ce^{-\gamma cn_1}=0.\]
Combining the three limit results for $S_1$, $S_2$ and $S_3$, it is easy to deduce that $\mathbf{P}$-a.s.,
\[\lim_{n\rightarrow\infty}|e^{-\gamma cn}X^c(cn)-m^c(\infty)|=0.\]

We continue by showing how the process $X$, when suitably scaled, converges along the subsequence $(cn)_{n\geq 0}$. Applying the conclusion of the previous paragraph, we find that $\mathbf{P}$-a.s.,
\begin{eqnarray}
\lefteqn{\limsup_{n\rightarrow\infty}|e^{-\gamma cn}X(cn)-m(\infty)|\leq}\nonumber\\
&& |m(\infty)-m^c(\infty)|+\limsup_{n\rightarrow\infty}e^{-\gamma cn}|X(cn)-X^c(cn)|.\label{upper}
\end{eqnarray}
Recall that the process $X^c$ and its discounted mean process $m^c$ depend on the integer $n_0$. By applying the dominated convergence theorem, it is straightforward to check that, as we let $n_0\rightarrow \infty$, the first of the terms in the above estimate, which is deterministic, converges to zero. For the second term, observe that
\begin{eqnarray*}
|X(t)-X^c(t)|&=&\sum_{i\in\Sigma_*}\eta_i(t+3\ln l(i))\mathbf{1}_{\{t+3\ln l(i)>cn_0 \}}\\
&\leq& 6 \#\{i\in\Sigma_*:\:t+ 3\ln l(i) > cn_0\}.
\end{eqnarray*}
Applying standard branching process results to the process described in the proof of Lemma \ref{props}, we are able to deduce the existence of a finite constant $C$ such that, as $t\rightarrow\infty$, we have $e^{-2t}\#\{i\in\Sigma_*:\:-\ln l(i) < t\}\rightarrow C$, $\mathbf{P}$-a.s., from which it follows that $\mathbf{P}$-a.s.,
\[\limsup_{n\rightarrow\infty}e^{-\gamma cn}|X(cn)-X^c(cn)|\leq 6C e^{-\gamma cn_0}.\]
Consequently, by choosing $n_0$ suitably large, the upper bound in (\ref{upper}) can be made arbitrarily small, which has as a result that $e^{-\gamma cn}X(cn)\rightarrow m(\infty)$ as $n\rightarrow \infty$, $\mathbf{P}$-a.s., for each $c$. The proposition is readily deduced from this using the monotonicity of $X$.
\end{proof}

\def\cprime{$'$}

\end{document}